\documentclass[12pt, twoside]{article}
\usepackage{amsmath,amsthm,amssymb}
\usepackage{times}
\usepackage{enumerate}

\theoremstyle{definition}

\usepackage{color}

\numberwithin{equation}{section}
\usepackage{tabular,array}
\usepackage{booktabs}
\usepackage{graphicx}
\usepackage{xspace}

\newcommand{\beq}{\begin{eqnarray*}}
\newcommand{\eeq}{\end{eqnarray*}}
\newcommand{\beqn}{\begin{eqnarray}}
\newcommand{\eeqn}{\end{eqnarray}}

\newcommand{\notinclude}[1]{}

\DeclareMathOperator{\PePoRuargmin}{\mathrm{arg \, min}}
\newcommand{\PePoRuvectz}[2]{\ensuremath{\left(\begin{array}{c}#1\\#2\end{array}\right)}}
\allowdisplaybreaks
\begin{document}
\baselineskip=17pt
\title{A Nested Variational Time Discretization for Parametric Anisotropic Willmore Flow}
\author{Ricardo Perl\\ Paola Pozzi\\
Martin Rumpf}
\date{}
\maketitle
\abstract{A variational time discretization of anisotropic Willmore flow combined with
a spatial discretization via piecewise affine finite elements is presented. 
Here, both the energy and the metric underlying the gradient flow are anisotropic, which in particular ensures that Wulff shapes are invariant up to scaling 
under the gradient flow.
In each time step of the gradient flow a 
nested optimization problem has to be solved.
Thereby, an outer variational problem reflects the time discretization of the
actual Willmore flow and involves an approximate anisotropic $L^2$-distance 
between two consecutive time steps and a fully implicit approximation of the anisotropic
Willmore energy. The anisotropic mean curvature needed to evaluate the energy integrand is replaced by the time discrete, approximate speed from an inner, fully implicit variational scheme for anisotropic mean curvature motion. 
To solve the nested optimization problem 
a Newton method for the associated Lagrangian is applied.
Computational results for the evolution of curves underline the robustness of the new scheme, in
particular with respect to large time steps.
}
\section{Introduction}
This paper generalizes a recently proposed variational time discretization \cite{OlRu10}  for 
isotropic Willmore flow to the corresponding anisotropic flow. Thereby, the anisotropic Willmore flow is defined as the gradient flow of the anisotropic Willmore energy
with respect to the corresponding anisotropic $L^2$-metric. 

The isotropic Willmore energy is given by 
$w[x] = \frac12 \int_{\mathcal M} {\mathbf h}^2 d\mathit{a}$, where $x$ denotes the identity map and ${\mathbf h}$ the mean curvature on a surface ${\mathcal M}$. The isotropic $L^2$-metric
is given by $(v,v)_{\mathcal M}= \int_{\mathcal{M}}|v|^2 d\mathit{a}$ , which is considered as a squared $L^2$-distance of the surface ${\mathcal M}$ being displaced with the vector field $v$ from the non displaced surface ${\mathcal M}$.
In the hypersurface case Willmore flow leads to a fourth order
parabolic evolution problem, which defines for a
given initial surface ${\mathcal M}_0$ a family of surfaces ${\mathcal M}(t)$ for $t\geq
0$ with ${\mathcal M}(0)={\mathcal M}_0$ \cite{Wi93,Si01,KuSc01}.  Applications of
a minimization of the isotropic Willmore energy and the corresponding Willmore flow 
include the processing of edge sets in imaging
\cite{NiMuSh93,Mum94,YoBe02,ChKaSh02}, geometry processing
\cite{WeWi92,BeSaCaBa00,BeBeSa01,PaXu06} and the mathematical treatment of biological
membranes \cite{He73,Se97,DuLiWa06}. Starting with work by Polden~\cite{Po95,Po96}
existence and regularity of Willmore flow was advanced in the last decade \cite{MaSi02, KuSc02c, Ri08}. 

Now, in the context of Finsler geometry the classical area functional is replaced by the 
anisotropic area functional ${\mathbf{a}_\gamma}[x]= \int_{\mathcal M}\gamma(n) d\mathit{a}$ with 
a local area weight $\gamma(n)$ depending on the local surface orientation. Here, $\gamma$ is 
a positive, $1$--homogeneous anisotropy function. 
In analogy to the isotropic case the anisotropic mean curvature
${\mathbf h}_\gamma$ is defined as the $L^2$--representation of the variation of the anisotropic area
in the direction of normal variations of the surface and can be evaluated as ${\mathbf h}_\gamma =
\mathrm{div}_{\mathcal M}\left(\nabla\gamma(n) \right)$. Hence, a possible first choice for an anisotropic Willmore functional 
is given by $\frac{1}{2} \,\int_{{\mathcal M}} {\mathbf h}_\gamma^2~d\mathit{a}\,$. Clarenz~\cite{Cl04} has
shown that Wulff shapes are the only minimizers of this energy. Palmer~\cite{Pa09} 
studied variational problems involving anisotropic bending energies
for surfaces with and without boundaries. 
Unfortunately, this energy definition does not imply the scale invariance property of 
Wulff shapes known for round spheres under isotropic Willmore flow.
Indeed, any round sphere is a stationary point of the isotropic Willmore functional in $\mathbb{R}^3$.
In $\mathbb{R}^2$ a circle of radius $R_0$ evolves under isotropic Willmore flow according the the 
ordinary differential equation $\dot R = \frac12 R^{-3}$. 
The counterpart of a round sphere in the anisotropic context is the Wulff shape as the unit ball with respect to 
the norm associated with the dual $\gamma^\ast$ of the anisotropy $\gamma$. 
But there is no such scaling law for the evolution of Wulff shapes under the above anisotropic variant of Willmore flow.

To ensure full consistency with the Finsler geometry, one has to adapt 
both the anisotropic energy and the anisotropic metric 
as suggested in \cite{Po12} (see Section \ref{B8:sec:2STD}).
Indeed, we make use of the associated anisotropic metric 
$\int_{{\mathcal M}} \gamma^\ast\left(v\right) (\nabla \gamma^\ast)\left(v\right) \cdot v  \, \gamma(n) d\mathit{a}$ (here only defined for $v(x) \neq 0$ for all $x \in {\mathcal M}$, cf. Section \ref{B8:sec:3STD} for the general case),
acting on a motion field $v$ of the surface ${\mathcal M}$ with normal $n$.
Furthermore, we will use the anisotropic area weight to define 
the anisotropic Willmore energy, i.e. $w_\gamma(x) = \frac{1}{2} \,\int_{{\mathcal M}} {\mathbf h}_\gamma^2\, \gamma(n) d\mathit{a}$.
Then, it actually turned out that Wulff shapes in $\mathbb{R}^2$ actually evolve according to the same evolution law for radial parameter valid for the evolution of circles under the isotropic flow.
Recently, Bellettini~\&~Mugnai~\cite{BeMu09} investigated the first variation of this functional in the smooth case.
Concerning the proper time and space discretization, this consistent choice of the anisotropic Willmore energy and the anisotropic metric on surface variations perfectly fits to the framework of the natural variational time discretization of geometric gradient flows. 

The finite element approximation of Willmore flow was first investigated by 
Rusu~\cite{Ru01} based on a mixed method for the 
surface parametrization $x$ and the mean curvature vector ${\mathbf h}\,n$ as
independent variables, see also \cite{ClDiDzRuRu04} for the 
application to surface restoration. In \cite{DrRu04b} a level set formulation of Willmore flow 
was proposed.
In the case of graph surfaces Deckelnick and
Dziuk \cite{DeDz06} were able to prove convergence of a related space discrete and time continuous
scheme. Deckelnick and Schieweck established convergence of a conforming finite element
approximation for axial symmetric surfaces \cite{DeSc09}. In the case of the elastic flow of curves
an error analysis was given by Dziuk and Deckelnick in \cite{DeDz09}.
An alternative scheme, which in particular ensures a better distribution of nodes on the evolving surface was 
presented by Barrett, Garcke and N\"urnberg \cite{BaGaNu07,BaGaN12}. Using discrete geometry calculus 
Bobenko and Schr\"oder \cite{BoSc05a} suggested a discrete Willmore flow of triangular surfaces. 
The time discretization of the second order, anisotropic mean curvature flow has been considered
by Dziuk already in \cite{Dz94,Dz99} and he gave convergence results for curves.
Diewald~\cite{Di05} has extended the discretization approach for isotropic Willmore flow 
of Rusu~\cite{Ru01}  to some anisotropic variant, for which Droske~\cite{Dr05} and Nemitz~\cite{Ne08}
investigated a level set discretization. 

Most of the above discretization methods are based on some semi-implicit time discretization, which requires 
the solution of linear systems of equations at each time step. Thereby, the involved geometric differential operators are
assembled on the surface from the previous time step. 
In the application one observes  strong restrictions on the time step size.
This shortcoming motivated the development of a new approach for the
time discretization of Willmore flow in \cite{OlRu10} based on the following 
general concept for a variational time discretization of gradient flows:
The gradient flow on a (in general infinite dimensional) manifold with respect to an energy $e[\cdot]$  and a metric $g$ on  the manifold,
is defined as the evolution problem $\dot x = - \mathrm{grad}_g e[x]$ with initial data
$x^0$, where $\mathrm{grad}_g e[x]$ is the representation of the variation $e'[x]$ in the metric $g$, i.e. 
$g(\mathrm{grad}_g e[x], \zeta) = e'[x](\zeta)$ for all infinitesimal variations $\zeta$ of $x$.
Now, one defines a time discrete family $(x^k)_{k=0,\cdots}$ with the desired property $x^k \approx x(k\tau)$ 
for the given time step size $\tau$. To this end, one successively solves a sequence of variational problems,
i.e. in time step $k$ 
\begin{eqnarray*} x^{k+1} &=& \PePoRuargmin\limits_{x} \mathrm{dist}(x^k,x)^2
+ 2\tau \, e[x]\,, 
\end{eqnarray*} 
where $\mathrm{dist}(x^k,x) =
\inf\limits_{\gamma \in \Gamma[x^k,x]} \int_0^1 \sqrt{g_{\gamma(s)}(\dot
\gamma(s),\dot \gamma(s))}d\mathit{s}$ denotes the Riemannian distance of $x$ from $x^k$ 
on the manifold and $\Gamma[x^k,x]$ is the set of smooth curves $\gamma$
with $\gamma(0)=x^k$ and $\gamma(1)=x$. 
The striking observation for this abstract scheme is that one immediately obtains an 
energy estimate, i.e. $e[x^{k+1}] + \frac1{2\tau} \mathrm{dist}(x^{k},x^{k+1})^2 \leq  e[x^k]\,.$ 
In the context of geometric flows, this approach was studied by Luckhaus and
Sturzenhecker \cite{LuSt95} leading to a fully implicit variational time discretization for mean curvature motion 
in $BV$ and by Chambolle~\cite{Ch04a}, who reformulated this scheme in terms of a level set method and generalized it for the approximation of anisotropic mean curvature motion 
in~\cite{BeCaChNo04,ChNo06}.
The time discretization for Willmore flow proposed in \cite{OlRu10} builds 
upon this general paradigm. 
In this paper, we will show how to adapt the approach to the time discretization of 
the anisotropic Willmore flow which is fully consistent with Finsler geometry.

The paper is organized as follows. In Section \ref{B8:sec:2STD} we briefly review the time discretization 
of isotropic Willmore flow. Building on these prerequisites 
the generalization to anisotropic Willmore flow is discussed in 
Section \ref{B8:sec:3STD}. Then, in Section \ref{B8:sec:space} we discuss 
a fully discrete numerical scheme based on piecewise affine
finite elements on simplicial surface meshes. In Section
\ref{B8:sec:opt} the Lagrangian calculus from PDE constraint optimization
is used to develop a suitable algorithm for the
solution of the nested optimization problem to be solved in each time step. 
Finally, in
Section \ref{B8:sec:numerics} computational results are presented. 
An appendix collects essential ingredients of the corresponding
algorithm.

\section{Review of the time discretization of isotropic Willmore flow}
\label{B8:sec:2STD} In this section we will briefly recall the nested time discretization of
isotropic Willmore from \cite{OlRu10}.
We denote a hypersurface in $\mathbb{R}^{d+1}$ by ${\mathcal M}={\mathcal M}[y]$. Here, 
$y$ indicates a parametrization of ${\mathcal M}$ and can also be considered 
as the identity map on ${\mathcal M}$ parametrizing ${\mathcal M}$ over itself.
Then, the abstract variational time discretization of
isotropic Willmore flow reads as follows: \\
For given surface ${\mathcal M}[x^k]$ with parametrization $x^k$ 
and a time step $\tau$ find a mapping $x=x[x^k]$ such that
$\mathrm{dist}({\mathcal M}[x^k],{\mathcal M}[x])^2 + \tau \int_{{\mathcal M}[x]} \mathbf{h}^2 d\mathit{a} \longrightarrow \min\,$, where 
\\$\mathrm{dist}({\mathcal M}[z],{\mathcal M}[v])^2 = \int_{{\mathcal M}[z]} (v-z)^2 d\mathit{a}$ is the squared
$L^2$-distance of surfaces ${\mathcal M}[v]$ from the surface ${\mathcal M}[z]$, $\mathbf{h}=\mathbf{h}[x]$
is the mean curvature of ${\mathcal M}[x]$, and
$\int_{{\mathcal M}[x]} d\mathit{a}$ denotes the surface area of ${\mathcal M}[x]$. 
Now, we take into account that the mean
curvature $\mathbf{h}=\mathbf{h}[x]$ is the $L^2$-gradient of the
area functional on a surface ${\mathcal M}[x]$ and that mean curvature motion is
the corresponding gradient flow. Thus, the mean curvature vector
$\mathbf{h}[x] n[x]$ with $n=n[x]$ denoting the normal on ${\mathcal M}[x]$
can be approximated by the discrete time derivative 
$\frac{y[x]-x}{\tilde \tau}$, where $y[x]$ is a suitable approximation of a single time step of 
the evolution of mean curvature motion with initial data $x$ and time step size $\tilde \tau$.
This time step itself can again be approximated using an (inner) variational scheme, i.e. we define $y[x]$ to be the minimizer of 
\begin{eqnarray}
\label{B8:eq:linMCM} {e_{\mbox{\tiny in}}}[x,y]:= \int_{{\mathcal M}[x]} (y-x)^2 +
\tilde\tau  |\nabla_{{\mathcal M}[x]} y|^2 d\mathit{a} \,. 
\end{eqnarray}
In fact, the corresponding Euler Lagrange equation is identical to the defining equation of 
the semi-implicit scheme for mean curvature motion proposed by Dziuk \cite{Dz91}:
\begin{eqnarray} \label{B8:eq:Dz}
0 = \int_{{\mathcal M}[x]} (y-x) \theta +
\tilde\tau  \nabla_{{\mathcal M}[x]} y \cdot \nabla_{{\mathcal M}[x]} \theta d\mathit{a}
\end{eqnarray}
Now, given $y[x]$ as the minimizer of \eqref{B8:eq:linMCM} for small $\tilde \tau$ the functional
$\frac12 \int_{{\mathcal M}[x]} \frac{(y[x]-x)^2}{\tilde \tau^2}d\mathit{a}$ is 
an approximation of the Willmore functional on ${\mathcal M}[x]$.
This approximation is then used to define a variational scheme for a time step of the actual Willmore flow. To this end, we consider for given surface parametrization $x^k$ the functional 
\begin{eqnarray*} {e_{\mbox{\tiny out}}}[x^k,x,y] :=
\int_{{\mathcal M}[x^k]} (x-x^k)^2 d\mathit{a} + \frac{\tau}{\tilde \tau^2}
\int_{{\mathcal M}[x]} (y-x)^2 d\mathit{a} \,,
\end{eqnarray*} 
where we suppose $y=y[x]$ to be the minimizer of \eqref{B8:eq:linMCM}.
To summarize, we obtain the following scheme for the $k$th time step of Willmore flow:\medskip

\emph{Given an initial surface ${\mathcal M}[x^0]$ with parametrization $x^0$ we
define a sequence of surfaces ${\mathcal M}[x^k]$ with parametrizations
$x^{k}$ for $k=1, \ldots$ via the solution of the following sequence of nested
variational problems
\begin{eqnarray}
x^{k+1} &=& \PePoRuargmin\limits_{x} {e_{\mbox{\tiny out}}}[x^k,x,y[x]] \, \mbox{, where}\label{B8:eq:outer} \\
y[x] &=& \PePoRuargmin\limits_{y} {e_{\mbox{\tiny in}}}[x,y] \label{B8:eq:inner} \,. \end{eqnarray}}
\vspace{2ex}

The inner variational problem \eqref{B8:eq:inner} is
quadratic, thus the resulting Euler--Lagrange equation \eqref{B8:eq:Dz} is  
linear and we end up with a PDE constrained
optimization problem to be solved in each time step. For more details we refer to 
\cite{OlRu10}.

\section{Nested time discretization for anisotropic Willmore flow}
\label{B8:sec:3STD}
Now, let us investigate the time discretization of anisotropic
Willmore flow in the co-dimension one case. Here, we will in particular focus on the proper choice of energy and metric.
We assume that $\gamma: \mathbb{R}^{d+1} \to [0,\infty)$ is a positive, 
$1$--homogeneous (i.e. $\gamma(\lambda p ) = |\lambda| \gamma(p)$ for all $\lambda \in \mathbb{R}, p \in \mathbb{R}^{d+1}$) and sufficiently regular function, that satisfies the ellipticity condition
\begin{equation}
 \gamma'' (p)qq \geq c_0 \| q\|^2 \quad \forall \, p,q \in \mathbb{R}^{d+1}, \| p \|=1, p \cdot q =0
\end{equation}
for some positive constant $c_0$ and the Euclidean norm $\| \cdot \|$. 
As already mentioned $\gamma(n)$ represents the anisotropic area weight for a surface normal $n$.
The isotropic case is recovered by choosing $\gamma(\cdot)= \| \cdot \|$.
We define the dual function of $\gamma$ as
\begin{eqnarray*}
\gamma^*(x) := \sup\{ \langle x,\psi \rangle ~|~ \psi \in B_{\gamma} \} \quad \forall ~ x \in \mathbb{R}^{d+1}\,,
\end{eqnarray*}
where $B_{\gamma}$ denotes the unit Ball in the $\gamma$-norm.
The ellipticity assumption ensures that $(\mathbb{R}^{d+1} ,\gamma)$
and its dual space $(\mathbb{R}^{d+1}, \gamma^*)$ are uniformly convex Banach spaces and the duality map $T: (\mathbb{R}^{d+1}, \gamma^*) \to (\mathbb{R}^{d+1}, \gamma)$, with
\begin{eqnarray*}
T(x)=\frac{1}{2} \partial (\gamma^*(x)^2),
\end{eqnarray*}
is an odd single-valued bijective continuous map. 
More precisely $T(0)=0$, $T(x)=\gamma^*(x)  \nabla \gamma^* (x)$ for $x \neq 0$, and $T^{-1}(\xi)= \gamma(\xi) \nabla \gamma (\xi)$ for $\xi \neq 0$.
For details we refer to \cite{Po12}. 
The unit ball $\mathcal{F}:=\{x \in \mathbb{R}^{d+1} \, : \,\gamma(x) \leq 1 \}$ 
in $(\mathbb{R}^{d+1} ,\gamma)$ is denoted the Frank diagram, the associated 
dual unit ball $\mathcal{W}:=\{x \in \mathbb{R}^{d+1} \, : \,\gamma^*(x) \leq 1 \}$ is the corresponding Wulff shape.
Wulff shapes are known to be solutions to the isoperimetric problem, that is 
$\partial \mathcal{W}$ minimizes the anisotropic area functional
\begin{equation}
 {\mathbf{a}_\gamma}[x]=\int_{{\mathcal M}[x]} \gamma(n[x]) d\mathit{a}
\end{equation}
(with $\gamma(n[x]) d\mathit{a}$ 
denoting the anisotropic area element) in the class of surfaces enclosing the same volume (cf. \cite{ClDzRu03} and the references therein). Now, based on the anisotropy $\gamma$ and its dual $\gamma^\ast$ we define an 
anisotropic distance  $\mathrm{dist}_\gamma$ of a manifold ${\mathcal M}[y]$ from a manifold ${\mathcal M}[x]$ 
by  
\begin{equation}\label{B8:eq:anisodist}
\mathrm{dist}_\gamma({\mathcal M}[x], {\mathcal M}[y])^2 := \int_{{\mathcal M}[x]} \gamma^*(y-x)^2  \gamma(n[x])d\mathit{a}
\end{equation}
for sufficiently regular $x$ and $y$. 
The choice of the norm $\gamma^*$ together with the anisotropic area weight $\gamma(n[x])$ in \eqref{B8:eq:anisodist}
reflects the fact that the anisotropic area  
of the boundary of a convex body $K\subset \mathbb{R}^{d+1}$ can be interpreted as
$${\mathbf{a}_\gamma} (\partial K)= \lim_{\epsilon \to 0} \frac{|K+\epsilon \mathcal{W}|-|K|}{\epsilon}\,,$$ 
where $|\cdot|$ denotes the usual Lebesgue volume in $\mathbb{R}^{d+1}$. 
In particular, the underlying 
metric structure is dictated by the Wulff shape and its norm $\gamma^*$ (see \cite{Be04}, \cite{Po12}  and references therein). 

Based on these considerations let us first consider anisotropic mean curvature motion, which is defined as the gradient flow of the anisotropic surface area with respect to the above anisotropic metric. In this case the variational time discretization 
is associated with the minimization of 
\begin{gather}\label{B8:eq:naturalanisocurvmotion}
\mathrm{dist}_\gamma({\mathcal M}[x], {\mathcal M}[y])^2 + 2 \tilde \tau \int_{{\mathcal M}[y]}
\gamma(n[y])d\mathit{a}\,
\end{gather}
with respect to $y$ for a given surface ${\mathcal M}[x]$ and $\tilde\tau > 0$. Let us denote by $y[x]$ the minimizer for given surface parameterization $x$.
The Euler Lagrange equation for \eqref{B8:eq:naturalanisocurvmotion} is given by
\begin{eqnarray} \nonumber
0 &=& \int_{{\mathcal M}[x]}  T(y-x) \cdot \theta  \, \gamma(n[x]) d\mathit{a} + \tilde \tau \langle a_\gamma' [y], \theta \rangle \\
&=& \tilde\tau \int_{{\mathcal M}[x]} T\left(\frac{y-x}{\tilde\tau}\right) \cdot \theta  \, \gamma(n[x]) d\mathit{a} 
+ \tilde \tau \langle a_\gamma' [y], \theta \rangle
\label{B8:eq:timediscreteMCM}
\end{eqnarray}
for smooth test functions $\theta:{\mathcal M}[x] \to \mathbb{R}^{d+1}$. Together with $\partial_t y(k\tilde\tau) \approx \frac{y-x}{\tilde\tau}$ this reflects the weak formulation of anisotropic mean curvature motion  given by 
\begin{equation}
\label{B8:eq:weakanisocurv}
 \int_{\mathcal{M}[y]}   T(\partial_t y) \cdot \theta \, \gamma(n[y]) d\mathit{a}= - \langle {\mathbf{a}_\gamma} ' [y], \theta \rangle
\end{equation}
for a parametrization $y$ and smooth test functions $\theta$ defined on ${\mathcal M}[y]$ (cf. \cite{Po12}).
Here, the variation of the anisotropic area functional is given by
$$\langle {\mathbf{a}_\gamma}'[y], \theta \rangle
=\int_{{\mathcal M}[y]} \mathbf{h}_\gamma[y] 
\frac{n[y]}{\gamma(n[y])} \cdot \theta \, \gamma(n[y]) d\mathit{a}\,,
$$ 
where  $\mathbf{h}_\gamma[y]= \mathrm{div}_{\mathcal M[y]}(n_\gamma[y])= \mathrm{div}_{\mathcal M[y]}(\nabla \gamma( n[y]))$ denotes the anisotropic mean curvature with $n_\gamma[y] = \nabla \gamma( n[y])$ (see \cite{Cl02}). 
Thus, from \eqref{B8:eq:weakanisocurv} we deduce that $T(\partial_t y) = -  \mathbf{h}_\gamma[y] \frac{n[y]}{\gamma(n[y])}$ or equivalently we achieve the strong formulation of  anisotropic mean curvature motion
$$\partial_t y =   \kappa_\gamma[y] := T^{-1}\left(- \mathbf{h}_\gamma[y] \frac{n[y]}{\gamma(n[y])}\right)
= -\mathbf{h}_\gamma[y] \nabla \gamma(n[y]) \,.$$ 
Indeed, as pointed out in \cite{Po12} the last equality holds due to the $1$--homogeneity of $\gamma$, i.e.
\begin{eqnarray*}
\gamma\left(- \mathbf{h}_\gamma[y] \frac{n[y]}{\gamma(n[y])}\right) \nabla\gamma\left(- \mathbf{h}_\gamma[y] \frac{n[y]}{\gamma(n[y])}\right)
&=& - \frac{\mathbf{h}_\gamma[y]}{\gamma(n[y])} \gamma(n[y]) \nabla \gamma(n[y]) \\
&=& - \mathbf{h}_\gamma[y] \nabla \gamma(n[y])\,.
\end{eqnarray*}

Next, we deal with the actual anisotropic Willmore flow and consider the anisotropic Willmore functional defined as follows for a  parametrization $x$ of ${\mathcal M}[x]$:
\begin{equation}
\label{B8:eq:Willmore}
 w_\gamma[x]:=\frac{1}{2} \int_{{\mathcal M}[x]} \mathbf{h}_\gamma[x]^2 \, \gamma(n[x]) d\mathit{a} 
 = \frac{1}{2} \int_{{\mathcal M}[x]} \gamma^* (\kappa_\gamma[x])^2 \, \gamma(n[x]) d\mathit{a} \,. 
\end{equation}
Here, we have used that the $1$--homogeneity and $\nabla \gamma(\xi) \in \partial \mathcal{W}$ for all $\xi \in \mathbb{R}^{d+1}$ imply 
$$
\gamma^* (\kappa_\gamma)^2 = \gamma^* \left(-\mathbf{h}_\gamma \nabla\gamma(n)\right)^2 = \mathbf{h}_\gamma^2 \gamma^*\left(- \nabla\gamma(n)\right)^2 = \mathbf{h}_\gamma^2\,.
$$
Then the abstract variational time discretization of anisotropic Willmore flow reads as follows: \medskip

Given ${\mathcal M}[x^k]$ and time step $\tau$ find a mapping $x=x[x^k]$ such that $x$ minimizes
\begin{gather}\label{B8:eq:naturalanisoWillmore}
\mathrm{dist}_\gamma({\mathcal M}[x^k], {\mathcal M}[x])^2 + \tau \int_{{\mathcal M}[x]} \gamma^* (\kappa_\gamma[x])^2 \, \gamma(n[x]) d\mathit{a}\,.
\end{gather}

As in the isotropic case, we will now replace the anisotropic mean curvature vector by the discrete speed 
extracted from a scheme for a single time step of anisotropic curvature flow \eqref{B8:eq:weakanisocurv}. In explicit, $\gamma^\ast(\frac{y[x]-x}{\tilde \tau})^2$ 
is a suitable approximation of $\mathbf{h}_\gamma^2[x] =\gamma^*(\kappa_\gamma[x])^2$, where $\frac{y[x]-x}{\tilde \tau}$ is the time discrete speed extracted from the
variational time discretization of anisotropic curvature motion.
Furthermore, we use the definition of the anisotropic distance measure in \eqref{B8:eq:anisodist}.
Finally, based on this approximation we derive the actual time discretization of anisotropic Willmore flow.
For a given surface parametrization $x^k$ of the surface ${\mathcal M}[x^k]$ at a time step $k$ 
we define the functionals
\begin{eqnarray*}
 {e_{\mbox{\tiny out}}}[x^k,x,y] &:=& \int_{{\mathcal M}[x^k]}  \gamma^*(x-x^k)^2 \, \gamma(n[x^k])d\mathit{a} + 
\frac{\tau}{\tilde{\tau}^2}  \int_{{\mathcal M}[x]} \gamma^*(y-x)^2 \, \gamma(n[x])d\mathit{a}\,, \\
{e_{\mbox{\tiny in}}} [x,y] &:=& \int_{{\mathcal M}[x]} \gamma^*(y-x)^2 \, \gamma(n[x])d\mathit{a} + 2\tilde \tau \int_{{\mathcal M}[y]} \gamma(n[y])d\mathit{a} \,,
\end{eqnarray*} 
and in analogy to the isotropic case above, we end up with the
following fully nonlinear variational time discretization of
anisotropic Willmore flow: \medskip

\emph{Given an initial surface ${\mathcal M}[x^0]$ with parametrization $x^0$ we
define a sequence of surfaces ${\mathcal M}[x^k]$ with parametrizations
$x^{k}$ for $k=1, \ldots$ via the solution of the following sequence of nested
variational problems
\begin{eqnarray}
\label{B8:eq:outer-ani}
x^{k+1} &=& \PePoRuargmin\limits_{x} {e_{\mbox{\tiny out}}}[x^k,x,y[x]], \,\quad \mbox{where} \\
\label{B8:eq:inner-ani}
y[x] &=&  \PePoRuargmin\limits_{y} {e_{\mbox{\tiny in}}} [x,y]  \,. 
\end{eqnarray}
}
Different from the variational scheme for isotropic Willmore flow, the inner variational problem is no longer quadratic.
It is worth to mention that this variational time discretization does not involve
derivatives of the anisotropy.  Nevertheless, as we will discuss below in the context of the actual computation,
differentiation is required to run Newton methods for the associated Lagrangian functional.
Indeed, for this we will need $\gamma, \gamma^* \in C^3(\mathbb{R}^{d+1} \setminus \{0 \})$; moreover, unless $(\gamma^*)^2 \in C^3(\mathbb{R}^{d+1})$ 
(which holds for $\gamma(p)= \sqrt{Ap \cdot p}$ with a symmetric positive definite matrix $A$), a regularization will be is required (see 
Section \ref{B8:sec:numerics} below).

Let us conclude this section with a study of boundaries $\partial \mathcal{W}$ of two-dimensional Wulff shapes  $\mathcal{W}$
moving under anisotropic Willmore flow in the plane.
To this end consider the parametrization $x: (0,T) \times S^1\to \mathbb{R}^2$, $x(t,\nu)=R(t)\nabla \gamma(\nu)$ of the boundary of a (rescaled) Wulff shape $R(t) \mathcal{W}$. Using the results given in 
\cite{Po12} it is easily seen that $x$ moves under anisotropic Willmore flow if $R(t)$ solves the ODE $$\dot{R(t)}=\frac{1}{2R(t)^3}\,.$$
Hence, we observe that Wulff shapes expand in time like in the isotropic case (cf. \cite{OlRu10}) with $R(t) = \sqrt[4]{R(0)^4 + 2t}$ . Next let us compare this with the time discrete evolution based on the proposed nested variational time discretization. 
We write $x,y,x^k:S^1 \to \mathbb{R}^2$, $x(\nu)=R \nabla \gamma(\nu)$, $y(\nu)=\tilde{R} \nabla \gamma(\nu)$, $x^k(\nu)= R^k \nabla \gamma(\nu)$.
Since $\gamma^*(\nabla \gamma(\nu))=1$ we immediately derive
\begin{eqnarray*}
 {e_{\mbox{\tiny out}}}[x^k,x,y] &=& (R-R^k)^2 {\mathbf{a}_\gamma}(x^k) +\frac{\tau}{\tilde{\tau}^2} (\tilde{R}-R)^2 {\mathbf{a}_\gamma}(x),\\
{e_{\mbox{\tiny in}}} [x,y] &=& (\tilde{R}-R)^2 {\mathbf{a}_\gamma}(x) +2\tilde{\tau} {\mathbf{a}_\gamma}(y).
\end{eqnarray*}
Considering variations  $y_\epsilon(\nu)=(\tilde{R} + \epsilon \psi) \nabla \gamma(\nu)$ in direction of the anisotropic normal $n_\gamma$ 
we infer from the inner problem that 
$$(\tilde{R}-R){\mathbf{a}_\gamma}(x)+\frac{\tilde{\tau}}{\tilde{R}}{\mathbf{a}_\gamma}(y)=0 .$$
More precisely, since ${\mathbf{a}_\gamma}(y)=\frac{\tilde{R}}{R} {\mathbf{a}_\gamma}(x)$ due to the homogeneity property of $\gamma$, we have that
$$ \tilde{R}=R-\frac{\tilde{\tau}}{R}.$$
This, together with ${\mathbf{a}_\gamma}(x)=\frac{R}{R^k} {\mathbf{a}_\gamma}(x^k)$, gives
$${e_{\mbox{\tiny out}}}[x^k,x,y]= {\mathbf{a}_\gamma}(x^k) \left( (R-R^k)^2 +\frac{\tau}{RR^k}\right),$$
from which we deduce
$$ \frac{R-R^k}{\tau}=\frac{1}{2R^k R^2}.$$
Note that this is a slightly different time step scheme than the one reported for the isotropic case ($\gamma(\cdot) = \|\cdot\|$)
in \cite[\S~2.1]{OlRu10}. This is due to fact that we use an implicit formulation of the inner problem as opposed to the linear equation \eqref{B8:eq:linMCM} in the scheme for isotropic Willmore flow (cf. Section \ref{B8:sec:2STD}).

\section{Finite element discretization in space}
\label{B8:sec:space}
Following the approach in \cite{OlRu10} we now derive a suitable spatial 
discretization based on piecewise affine finite elements. This is in close correspondence to the surface finite element approach by Dziuk
\cite{Dz88}. To this end we consider simplicial meshes ${\mathcal M}[X]$ 
as approximations of the hypersurfaces ${\mathcal M}[x]$ in $\mathbb{R}^{d+1}$, i.e. 
polygonal curves for $d=1$ and triangular surfaces for $d=2$. Thereby, $X$
is a parametrization of the simplicial mesh ${\mathcal M}[X]$ which is uniquely
described by a vector $\bar X$ of vertex positions of the mesh.  Here, and in what follows, 
we will always denote discrete quantities with
upper case letters to distinguish them from the corresponding
continuous quantities in lower case letters. Furthermore, a bar on
top of a discrete function indicates the associated vector of nodal values,
i.e. $\bar X = (\bar X_i)_{i\in I}$, where $\bar X_i = (X_i^1,\cdots, X_i^{d+1})$ 
is the coordinate vector of the $i$th vertex of the
mesh and $I$ denotes the index set of vertices. For $d=1$ each
element $T$ is a line segment with nodes $X_0$
and $X_1$ (using local indices) and for $d=2$ the elements $T$ 
are planar triangles with vertices $X_0$, $X_1$, and
$X_2$ and edge vectors $F_0=X_2-X_1$, $F_1=X_0-X_2$, and
$F_2=X_1-X_0$.
Given a simplicial surface ${\mathcal M}[X]$, the associated piecewise affine finite element space
is given by 
$${\mathcal V}({\mathcal M}[X]) := \left\{ U \in C^0({\mathcal M}[X])\, | \, U|_T \in {\mathcal P}_1 \, \forall T \in {\mathcal M}[X] \right\} $$
with the nodal basis denoted by $\{\Phi_i\}_{i\in I}$. Here, ${\mathcal P}_1$ is the space of affine functions on a simplex $T$.
Thus, for $U\in {\mathcal V}({\mathcal M}[X])$ we obtain $U=\sum_{i\in I} U(X_i) \Phi_i$ 
and $\bar U = (U(X_i))_{i\in I}$.
Let us emphasize, that the parametrization mapping $X$ itself is considered as an
element  in ${\mathcal V}({\mathcal M}[X])^{d+1}$ and we recover the vector of nodes $\bar X= (X_i)_{i\in I}$.

With these algorithmic ingredients at hand
we now can derive a fully discrete nested time discretization of
anisotropic Willmore flow, as  the spatially discrete counterpart of
\eqref{B8:eq:outer-ani} and \eqref{B8:eq:inner-ani}:

\emph{Given a discrete initial surface ${\mathcal M}[X^0]$ with discrete
parametrization $X^0$ we compute a sequence of surfaces ${\mathcal M}[X^k]$
with parametrizations $X^{k}$ by solving the nested, finite dimensional variational problems 
\begin{eqnarray}
X^{k+1} &=& \PePoRuargmin\limits_{X\in {\mathcal V}({\mathcal M}[X^k])^{d+1}} {\mathcal{E}_{\mbox{\tiny out}}}[X^k,X,Y[X]], \quad \mbox{where} \label{B8:eq:OUTER}\,,\\
Y[X] &=& \PePoRuargmin\limits_{Y\in {\mathcal V}({\mathcal M}[X])^{d+1}} {\mathcal{E}_{\mbox{\tiny in}}}[X,Y] \,.
\label{B8:eq:INNER} 
\end{eqnarray}
}
Here, the functionals ${\mathcal{E}_{\mbox{\tiny in}}}$ and ${\mathcal{E}_{\mbox{\tiny out}}}$ are straightforward spatially discrete counterpart
of the functionals ${e_{\mbox{\tiny in}}}[x,y]$ and ${e_{\mbox{\tiny out}}}[x^k,x,y]$ and defined by
\begin{eqnarray*}
{\mathcal{E}_{\mbox{\tiny in}}}[X,Y] &:=& \int_{{\mathcal M}[X]} \mathbf{I}\left(\gamma^*(Y-X)^2\right) \, \gamma(N[X]) d\mathit{a} 
+2 \tilde{\tau}  \int_{{\mathcal M}[Y]}  \gamma(N[Y]) d\mathit{a}\,,\\[2mm]
{\mathcal{E}_{\mbox{\tiny out}}}[X^k,X,Y] &:=& \int_{{\mathcal M}[X^k]} \!\!\mathbf{I}\left(\gamma^*(X-X^k)^2\right) \, \gamma(N[X^k])d\mathit{a} \\
&&+ \frac{\tau}{\tilde \tau^2} \int_{{\mathcal M}[X]} \mathbf{I}\left(\gamma^*(Y-X)^2\right) \,\gamma(N[X])d\mathit{a}\,,
\end{eqnarray*} 
where the nodal interpolation operator  $\mathbf{I}$  renders the resulting scheme fully practical.
To simplify the exposition, we introduce the discrete quadratic form ${\mathbf{M}_\gamma}[Z,X] = \int_{{\mathcal M}[X]} \mathbf{I}\left(\gamma^\ast(Z)^2\right) \gamma(N[X]) d\mathit{a}$ (a nonlinear counter part of the quadratic form induced by the lumped mass matrix) and the discrete anisotropic area functional ${\mathbf{A}_\gamma}[Y] = \int_{{\mathcal M}[Y]} \gamma(N[Y]) d\mathit{a}\,,$
both of which are assembled from local contributions on simplices of the underlying simplicial grid $\mathcal{T}_h$:
\begin{eqnarray}
{\mathbf{M}_\gamma}[Z,X] &=& \sum_{T \in \mathcal{T}_h} \frac{1}{(d+1)!} \left(\sum_{i =0,\ldots, d} \gamma^\ast(\bar Z_{T,i})^2 \right) \gamma(R_T[\bar X])
\label{B8:Mgamma} \\
{\mathbf{A}_\gamma}[X] &=&  \sum_{T \in \mathcal{T}_h} \frac1{d!} \gamma(R_T[\bar X])
\label{B8:Agamma}
\end{eqnarray}
Here, $R_T [\bar X] = D^{90} (\bar X_{T,1}- \bar X_{T,0})$ for $d=1$ and $R_T [\bar X] =  (\bar X_{T,1}- \bar X_{T,0}) \wedge (\bar X_{T,2}- \bar X_{T,0})$ for $d=2$.
Hence, we can rewrite 
\begin{eqnarray*}
{\mathcal{E}_{\mbox{\tiny out}}}[X^k,X,Y] &=& {\mathbf{M}_\gamma}[X-X^k,X^k] + \frac{\tau}{\tilde \tau^2} {\mathbf{M}_\gamma}[Y-X,X]\,, \\
{\mathcal{E}_{\mbox{\tiny in}}}[X,Y] &=& {\mathbf{M}_\gamma}[Y-X,X] + 2 \tilde \tau {\mathbf{A}_\gamma}[Y]\,.
\end{eqnarray*}
The  necessary condition for $Y[X]$ to be a minimizer of ${\mathcal{E}_{\mbox{\tiny in}}}[X,\cdot]$ is given by the corresponding discrete Euler 
Lagrange equation
\begin{eqnarray*}
0 = \partial_Y {\mathcal{E}_{\mbox{\tiny in}}}[X,Y[X]](\Theta) = \partial_Z {\mathbf{M}_\gamma}[Y-X,X](\Theta) + 2 \tilde \tau \partial_Y{\mathbf{A}_\gamma}[Y](\Theta)
\end{eqnarray*}
for all $\Theta \in {\mathcal V}({\mathcal M}[X])^{d+1}$.

\section{Optimization algorithm for the time steps}
\label{B8:sec:opt}
In this section, the actual optimization algorithm for the nested, fully discrete 
variational problem derived in Section \ref{B8:sec:space} is presented.
Thereby, we apply a step size controlled Newton method (cf. \cite{ScWe04} section 7) for the corresponding Lagrangian
(cf. Nocedal~\&~Wright~\cite{NoWr99}).
In our context the Lagrangian function for problem~\eqref{B8:eq:OUTER}, \eqref{B8:eq:INNER}
is given by
$$
  \mathcal{L}[\bar X,\bar Y,\bar P]={\mathcal{E}_{\mbox{\tiny out}}}[X^k,X,Y]-\partial_Y {\mathcal{E}_{\mbox{\tiny in}}}[X,Y](P)
$$
for independent unknowns $\bar X,\bar Y\in \mathbb{R}^{(d+1) |I|}$ and the
Lagrange multiplier $\bar P \in \mathbb{R}^{(d+1) |I|}$ (with a slight
misuse of notation, we consider these unknowns as finite element function
in the spaces ${\mathcal V}({\mathcal M}[X^k])^{d+1}$ and ${\mathcal V}({\mathcal M}[X])^{d+1}$, respectively, or as
the associated nodal vector in $\mathbb{R}^{(d+1) |I|}$). 
For an extensive discussion of the Lagrangian ansatz we refer to \cite{OlRu09}.
Now, we ask for critical points $(\bar X,\bar Y,\bar P)$ of $L$.
Indeed, $0=\partial_{\bar P} \mathcal{L}[\bar X,\bar Y,\bar P](\bar \Theta) =
- \partial_Y {\mathcal{E}_{\mbox{\tiny in}}}[X,Y](\Theta)$ is the Euler Lagrange equation of the
inner minimization problem with respect to $Y$ for given $X$ and $0= \partial_{\bar
Y} \mathcal{L}[\bar X,\bar Y,\bar P](\bar  \Theta) =\partial_{Y}{\mathcal{E}_{\mbox{\tiny out}}}[X^k,X,Y](\Theta) -
\partial_Y^2 {\mathcal{E}_{\mbox{\tiny in}}}[X,Y](P,\Theta)$ is the defining equation for the dual
solution $P$ given $Y$ as the solution of the above Euler Lagrange equation.
Finally, the Euler Lagrange equation for the actual constraint optimization problem 
coincides with $$0= \partial_{\bar X} \mathcal{L}[\bar X,\bar Y,\bar P](\bar  \Theta) =
\partial_X {\mathcal{E}_{\mbox{\tiny out}}}(X^k,X,Y)(\Theta) - \partial_X\partial_Y {\mathcal{E}_{\mbox{\tiny in}}}[X,Y](P,\Theta)\,.$$
For the gradient of the Lagrangian $\mathcal{L}$ we obtain
\begin{eqnarray*} 
\mathrm{grad}\, \mathcal{L} = \left(\begin{array}{c}
\partial_X {\mathcal{E}_{\mbox{\tiny out}}} - \partial_X \partial_Y {\mathcal{E}_{\mbox{\tiny in}}}(P) \\
\partial_Y {\mathcal{E}_{\mbox{\tiny out}}} - \partial_Y^2 {\mathcal{E}_{\mbox{\tiny in}}}(P)\\ 
-\partial_Y {\mathcal{E}_{\mbox{\tiny in}}}
\end{array}\right)
\end{eqnarray*} 

with
\begin{eqnarray*}
\partial_X {\mathcal{E}_{\mbox{\tiny out}}} [X^k,X,Y](\Theta) &=& 
\partial_Z {\mathbf{M}_\gamma}[X-X^k,X^k](\Theta) \nonumber \\
&& + \frac{\tau}{\tilde \tau^2} ( \partial_X {\mathbf{M}_\gamma}[Y-X,X](\Theta) - \partial_Z {\mathbf{M}_\gamma}[Y-X,X](\Theta))\,,
\\
\partial_Y {\mathcal{E}_{\mbox{\tiny out}}} [X^k,X,Y](\Theta) &=&  \frac{\tau}{\tilde \tau^2}  \partial_Z {\mathbf{M}_\gamma}[Y-X,X](\Theta) \,,
\\
\partial_X  \partial_Y {\mathcal{E}_{\mbox{\tiny in}}}[X,Y](P,\Theta) &=& 
-\partial_Z^2 {\mathbf{M}_\gamma}[Y-X,X](P,\Theta) + \partial_X \partial_Z {\mathbf{M}_\gamma}[Y-X,X](P,\Theta) \,,
\\
\partial_Y^2 {\mathcal{E}_{\mbox{\tiny in}}}[X,Y](P,\Theta) &=&
\partial^2_Z {\mathbf{M}_\gamma}[Y-X,X](P,\Theta) + 2 \tilde \tau \partial^2_Y{\mathbf{A}_\gamma}[Y](P,\Theta)\,.
\end{eqnarray*}
The Hessian of $\mathcal{L}$, which is required to implement a Newton scheme, is given
(in abbreviated form) by \begin{eqnarray*} \mathrm{Hess}\, \mathcal{L} = \left(\begin{array}{ccc}
\partial_X^2 {\mathcal{E}_{\mbox{\tiny out}}} - \partial_X^2 \partial_Y {\mathcal{E}_{\mbox{\tiny in}}}(P) &
\partial_X\partial_Y {\mathcal{E}_{\mbox{\tiny out}}} - \partial_X \partial_Y^2 {\mathcal{E}_{\mbox{\tiny in}}}(P) &
-\partial_X \partial_Y {\mathcal{E}_{\mbox{\tiny in}}} \\ 
\partial_X\partial_Y {\mathcal{E}_{\mbox{\tiny out}}} - \partial_X \partial_Y^2 {\mathcal{E}_{\mbox{\tiny in}}}(P) &
\partial_Y^2 {\mathcal{E}_{\mbox{\tiny out}}} - \partial_Y^3 {\mathcal{E}_{\mbox{\tiny in}}}(P) &
- \partial_Y^2 {\mathcal{E}_{\mbox{\tiny in}}} \\ 
-\partial_X \partial_Y {\mathcal{E}_{\mbox{\tiny in}}} &
- \partial_Y^2 {\mathcal{E}_{\mbox{\tiny in}}} &
0
\end{array}\right)\,.
\end{eqnarray*} 
The different terms in $\mathrm{Hess}\,  \mathcal{L}$ are evaluated as follows: 
\begin{eqnarray*}
\partial_X^2 {\mathcal{E}_{\mbox{\tiny out}}}(\Theta,\Psi) &=& 
\partial_Z^2 {\mathbf{M}_\gamma}[X-X^k,X^k](\Theta,\Psi) + \frac{\tau}{\tilde \tau^2} \big( \partial_X^2 {\mathbf{M}_\gamma}[Y-X,X](\Theta,\Psi) \\
& &-  2 \partial_Z \partial_X {\mathbf{M}_\gamma}[Y-X,X](\Theta,\Psi) + \partial_Z^2 {\mathbf{M}_\gamma}[Y-X,X](\Theta,\Psi) \big)
\,,\\
\partial_Y \partial_X  {\mathcal{E}_{\mbox{\tiny out}}}(\Theta,\Psi) &=& 
\frac{\tau}{\tilde \tau^2} \big( \partial_Z \partial_X  {\mathbf{M}_\gamma}[Y-X,X](\Theta,\Psi)
-  \partial_Z^2  {\mathbf{M}_\gamma}[Y-X,X](\Theta,\Psi) \big)
\,,\\
\partial_Y^2 {\mathcal{E}_{\mbox{\tiny out}}}(\Theta,\Psi) &=& 
\frac{\tau}{\tilde \tau^2}  \partial_Z^2 {\mathbf{M}_\gamma}[Y-X,X](\Theta,\Psi) 
\,,\\
\partial_X^2 \partial_Y {\mathcal{E}_{\mbox{\tiny in}}}(\Theta,\Psi,\Xi) &=&
\partial_Z^3 {\mathbf{M}_\gamma}[Y-X,X](\Theta,\Psi,\Xi) -\partial_X \partial_Z^2 {\mathbf{M}_\gamma}[Y-X,X](\Theta,\Psi,\Xi) \\
&&\!\!\! - \partial_X \partial_Z^2 {\mathbf{M}_\gamma}[Y-X,X](\Theta,\Xi,\Psi) +\partial_X^2 \partial_Z {\mathbf{M}_\gamma}[Y-X,X](\Theta,\Psi,\Xi)
\,,\\
\partial_X \partial_Y^2 {\mathcal{E}_{\mbox{\tiny in}}}(\Theta,\Psi,\Xi) &=& 
-\partial^3_Z {\mathbf{M}_\gamma}[Y-X,X](\Theta,\Psi,\Xi) 
+\partial_X \partial^2_Z {\mathbf{M}_\gamma}[Y-X,X](\Theta,\Psi,\Xi) 
\,,\\
\partial_Y^3 {\mathcal{E}_{\mbox{\tiny in}}}(\Theta,\Psi,\Xi) &=&  
\partial^3_Z {\mathbf{M}_\gamma}[Y-X,X](\Theta,\Psi,\Xi) + 2\tilde \tau \partial^3_Y{\mathbf{A}_\gamma}[Y](\Theta,\Psi,\Xi)
\,.\\
\end{eqnarray*}
In the implementation of the proposed scheme it is convenient to directly treat the squared, dual anisotropy $\gamma^{*,2}(.) := (\gamma^*(.))^2$
in the calculation of derivatives of the anisotropic functionals, which is particularly advantageous  for anisotropies of the type 
$\gamma(p) = \sum_{k=1}^K\sqrt{ p \cdot G_{k}p}$
where the $G_k$ are symmetric and positive definite (cf. Garcke et al. \cite{BaGaN08}). 
The different terms of the gradient $\mathrm{grad}\, \mathcal{L}$ and the Hessian $\mathrm{Hess}\, \mathcal{L}$ are in the usual way assembled from local contribution on simplices of the polygonal mesh. The required formulas  are given in the Appendix.

\section{Numerical results}
\label{B8:sec:numerics}
In this section, we show applications of the 
proposed algorithm to the evolution of curves in $\mathbb{R}^2$ under 
anisotropic Willmore flow. 
Beside anisotropies with ellipsoidal Wulff shapes we study regularized crystalline anisotropies $\gamma(\cdot) =  \|\cdot\|_{\ell^1}$ and 
$\gamma(\cdot) =  \|\cdot\|_{\ell_\infty}$ based on a suitable regularization. A particular emphasis is on the verification of the robustness and stability of the proposed approach in particular for large time steps. 
Furthermore, we experimentally verify that Wulff shapes grow self-similar in time under the corresponding anisotropic Willmore flow.
\begin{figure}[h]
\begin{tabular}{m{0.2\linewidth}m{0.3\linewidth}m{0.4\linewidth}}
\includegraphics[width=1\linewidth]{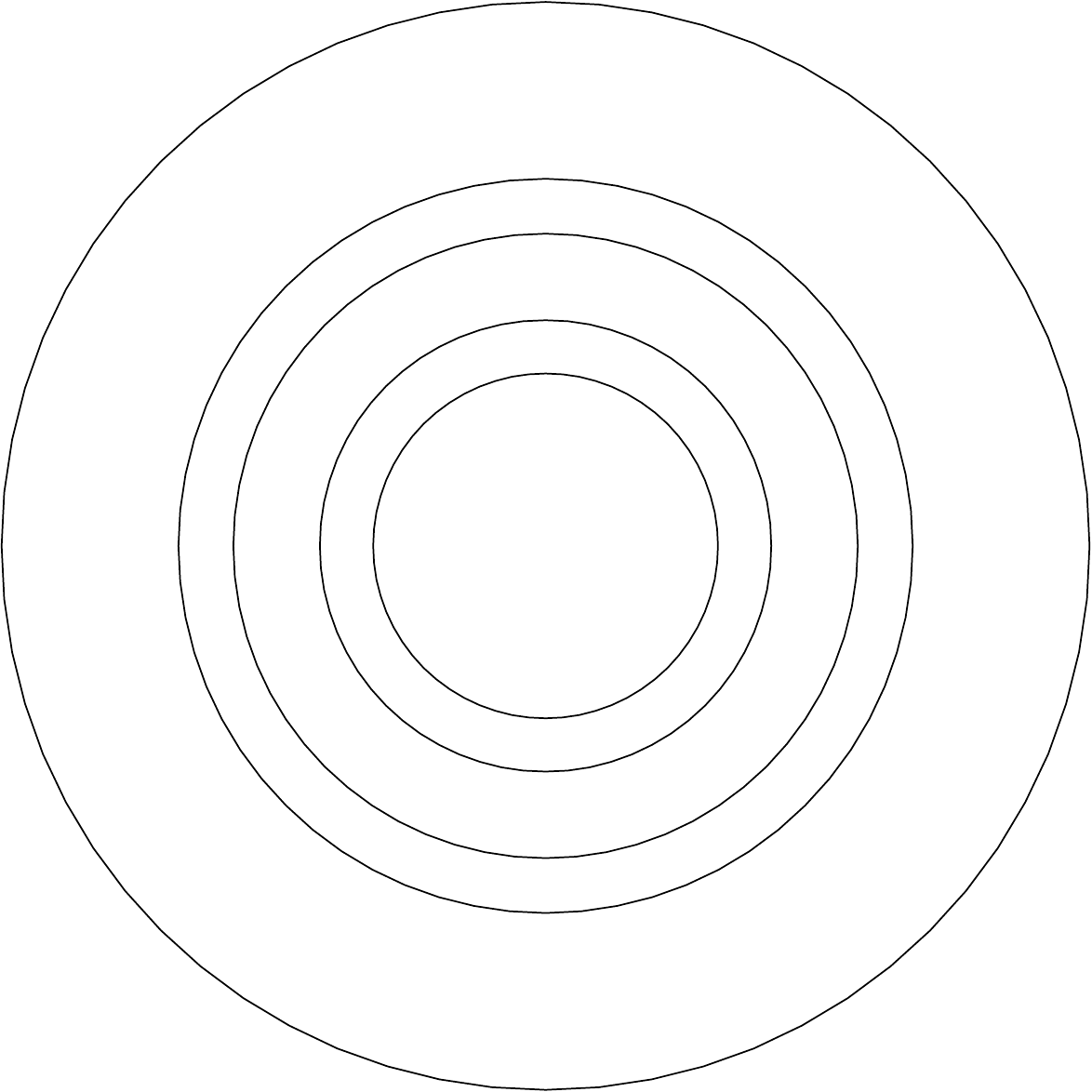} & \includegraphics[width=1\linewidth]{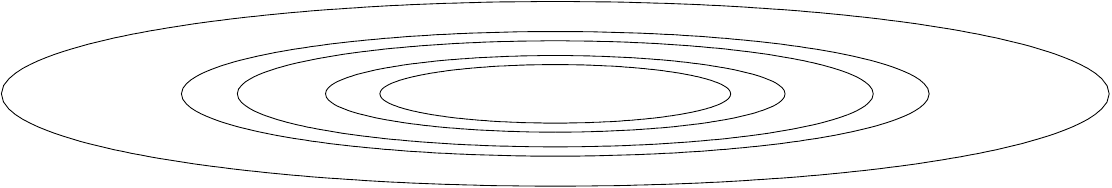} & \includegraphics[width=1\linewidth]{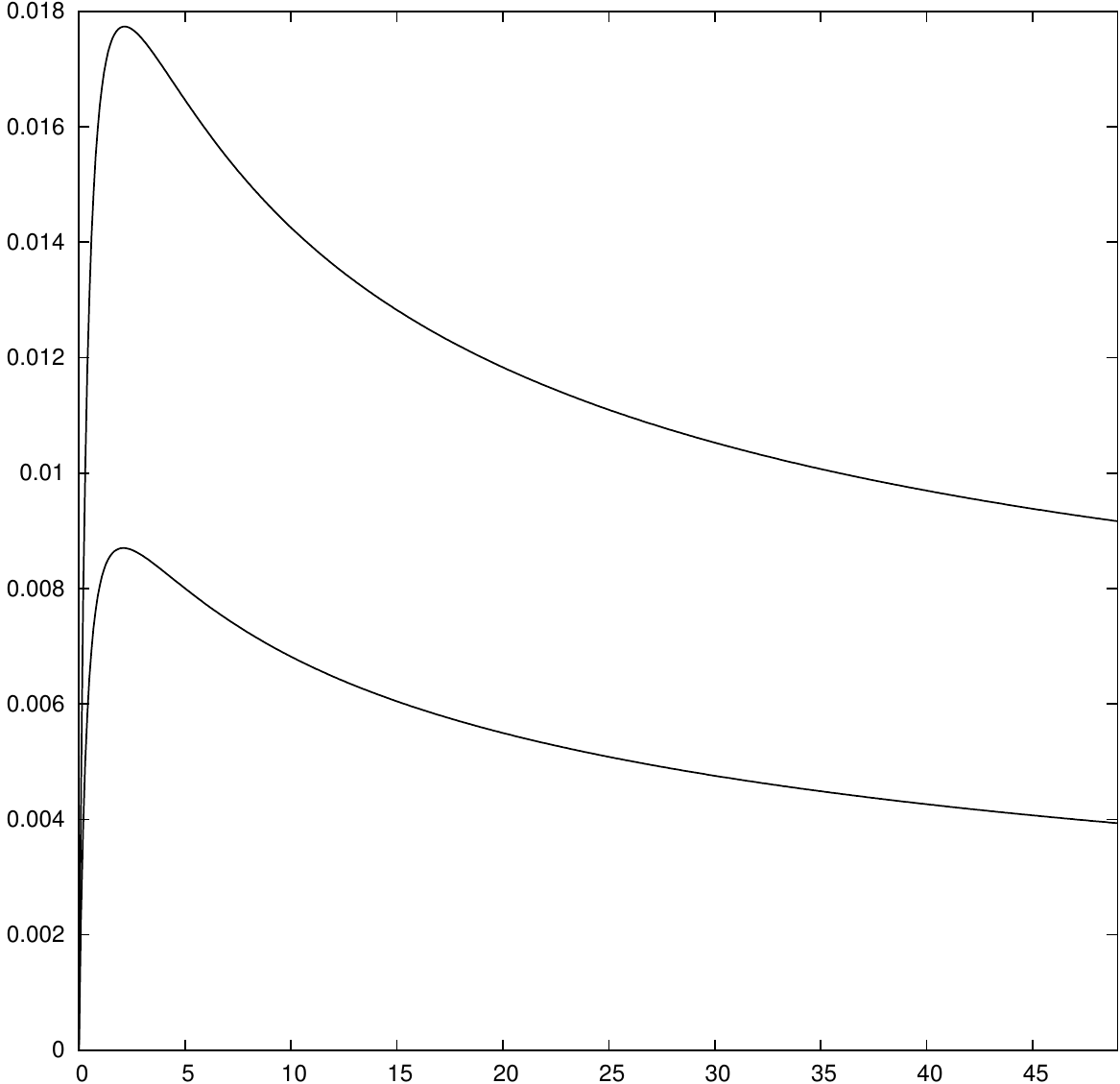}
\end{tabular}
\caption{The evolution of an unit circle under isotropic Willmore flow is plotted on the left. For the computation we used as initial grid size $h = 0.0981$ resulting from $64$ vertices. 
Furthermore, $\tau = h$, $\tilde\tau = h^2$ and the resulting discrete curves are shown for $t = 0, 10\tau, 50\tau,100\tau,500\tau$.
In the middle we display the evolution of an ellipse (with half axes 6 and 1) under anisotropic Willmore flow with $256$ elements and $h = 0.0984$. 
Here, we consider $\tau = h$, $\tilde\tau = h^2$ and display the approximate solutions for  $t = 0, 10\tau, 50\tau,100\tau,500\tau$.
Next, the associated $L^2$-errors are plotted over time on the right, where the lower error curve corresponds to the evolution results on the left.}
\label{B8:fig:ellipse_evolution}
\end{figure}

At first, we study anisotropies of the type
\begin{eqnarray*}
\gamma(z) = \sqrt{a_1^2 z_1^2 + a_2^2 z_2^2}
\end{eqnarray*}
for given $a_1,a_2 > 0$. In that case the squared dual anisotropy function is given by
\begin{eqnarray*}
\gamma^{*,2}(z) = \frac{z_1^2}{a_1^2} + \frac{z_2^2}{a_2^2}.
\end{eqnarray*}
Figure \ref{B8:fig:ellipse_evolution}  compares the evolution of a circle of radius $R_0=1$ under 
isotropic Willmore flow for $a_1 = a_2 = 1$ with the evolution of
an ellipse with half axes $a_1=6$ and $a_2=1$ under the corresponding anisotropic flow.
As discussed in Section \ref{B8:sec:2STD} in both cases the initial curve ${\mathcal M}_0$ expand in a self-similar fashion, i.e. $\mathcal{M}[x(t)] = R(t) {\mathcal M}_0$ with $R(t) = \sqrt[4]{R_0^4 + 2t}$ for $r_0 > 0$.
In Figure \ref{B8:fig:ellipse_evolution} we plot the evolution of the error 
$err(h) := \| \mathcal{I}_h x(t) - x_{h}(t) \|_{L^2}$ in time. Thereby, the $L^2$-error is evaluated on the 
polygonal curve $x_h(t)$ and  $\mathcal{I}_h$ denotes the nodal interpolation of $x(t)$ at the projected positions of the nodes of $x_h(t)$ in direction $\nabla \gamma(n[x_h(t)])$.
In Table \ref{B8:table:eoc1} and \ref{B8:table:eoc2} we provide results on the experimental order of convergence 
$\textit{eoc} := \log({err(h_1)}/{err(h_2)}) / \log({h_1}/{h_2})$ for varying grid and time step size in case of the evolution of the circle and the ellipse. 

\begin{table}[h]
\label{B8:tab:circle}
\begin{center}
\begin{tabular}{ | c || c | c | c || c | c | c |}
 \hline
  & & $L^2$-error & & & $L^2$-error & \\
  $n$   & $h(t)$ & $(\tau = \tilde{\tau} = h_0^2)$ & $eoc$ & $h(t)$ & $(\tau = \tilde{\tau} = h_0)$ & $eoc$  \\ \hline
  4   & 4.166e-1 & 4.830e-3    &       & 4.482e-1 & 1.916e-2 &  \\
  5   & 2.096e-1 & 1.328e-3    & 1.879 & 2.258e-1 & 1.087e-2 & 0.826 \\
  6   & 1.049e-1 & 3.403e-4    & 1.969 & 1.132e-1 & 5.804e-3 & 0.909 \\
  7   & 5.249e-2 & 8.561e-5    & 1.992 & 5.668e-2 & 3.000e-3 & 0.954 \\
  8   & 2.625e-2 & 2.144e-5    & 1.998 & 2.836e-2 & 1.525e-3 & 0.977 \\
\hline
\end{tabular}
\end{center}
\caption{\label{B8:table:eoc1}The $L^2$-error between the exact solution of the self-similar evolution of circles under Willmore flow and the discrete solution of the fully implicit variational time discretization is plotted at time $t = 0.1542$ for a grid size $h(t)$ (left) and $t = 0.3927$ (right). On the left we consider time step sizes $\tau$ and $\tilde \tau$ of the order of the squared spatial grid size $h_0$ at the initial time $0$, whereas on the right both time step sizes are taken equal to the grid size. In both cases we have considered $2^n$ vertices for the polygon, resulting in an initial grid size $h_0 = \frac{2\pi}{2^n}$.}
\end{table}

\begin{figure}[h]
\begin{center}
\includegraphics[width=0.3\linewidth]{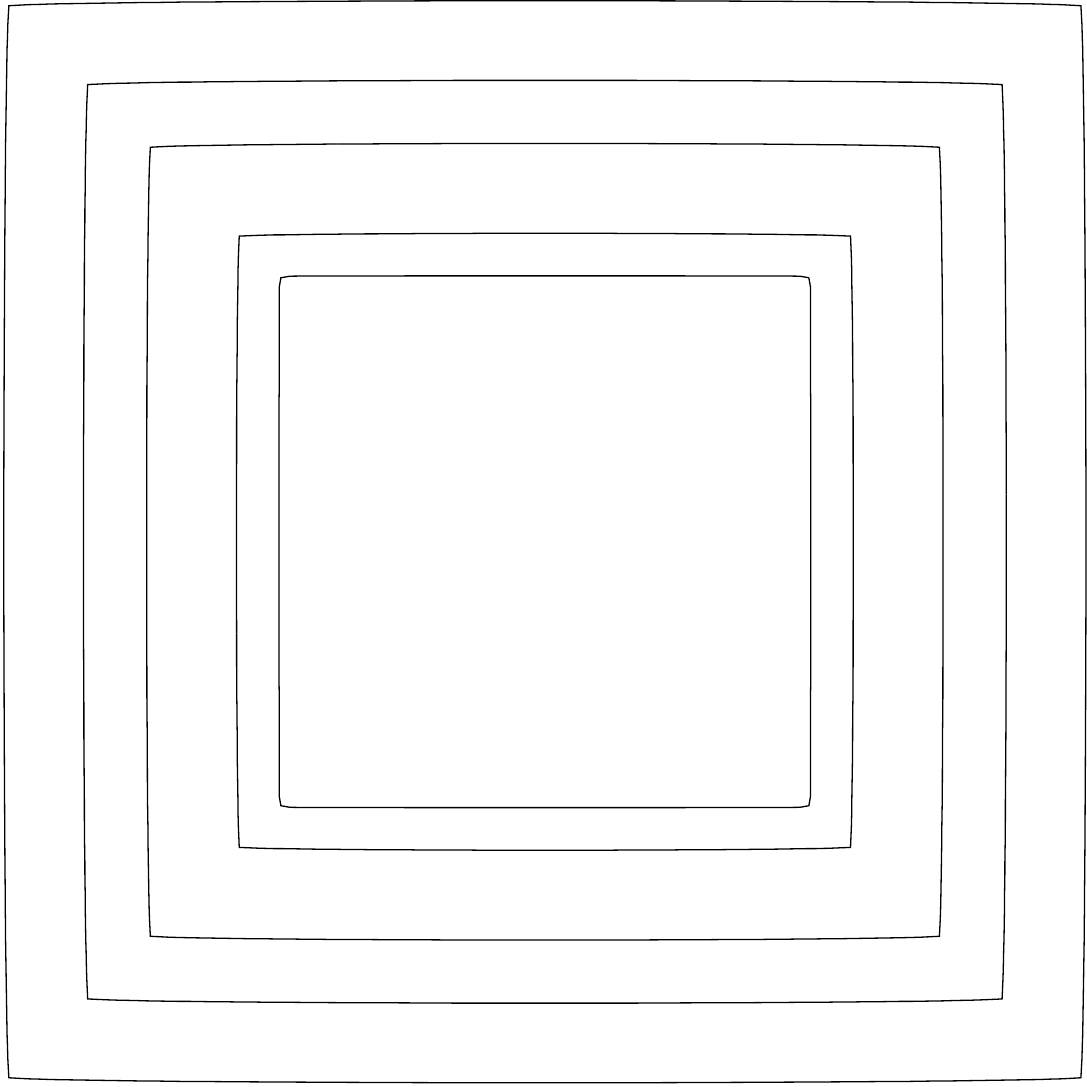}
\end{center}
\caption{ \label{B8:fig:evogammalinf} Evolution of the unit sphere with respect to the regularized $\ell^{\infty}$-norm under anisotropic Willmore flow for the anisotropy $\|\cdot\|_{\ell^1_{\varepsilon}}$ with $\varepsilon = 0.0001$. For this computation we consider $200$ vertices
leading to an initial grid size $h_0 = 0.04$. Furthermore, $\tau = h_0$ and $\tilde\tau = h_0^2$ and the resulting discrete curves are shown for $t = 0, 10\tau, 50\tau, 100\tau, 200\tau$.}
\end{figure}

\begin{figure}[h]
\begin{tabular}{m{0.25\linewidth}m{0.25\linewidth}m{0.4\linewidth}}
\includegraphics[width=1\linewidth]{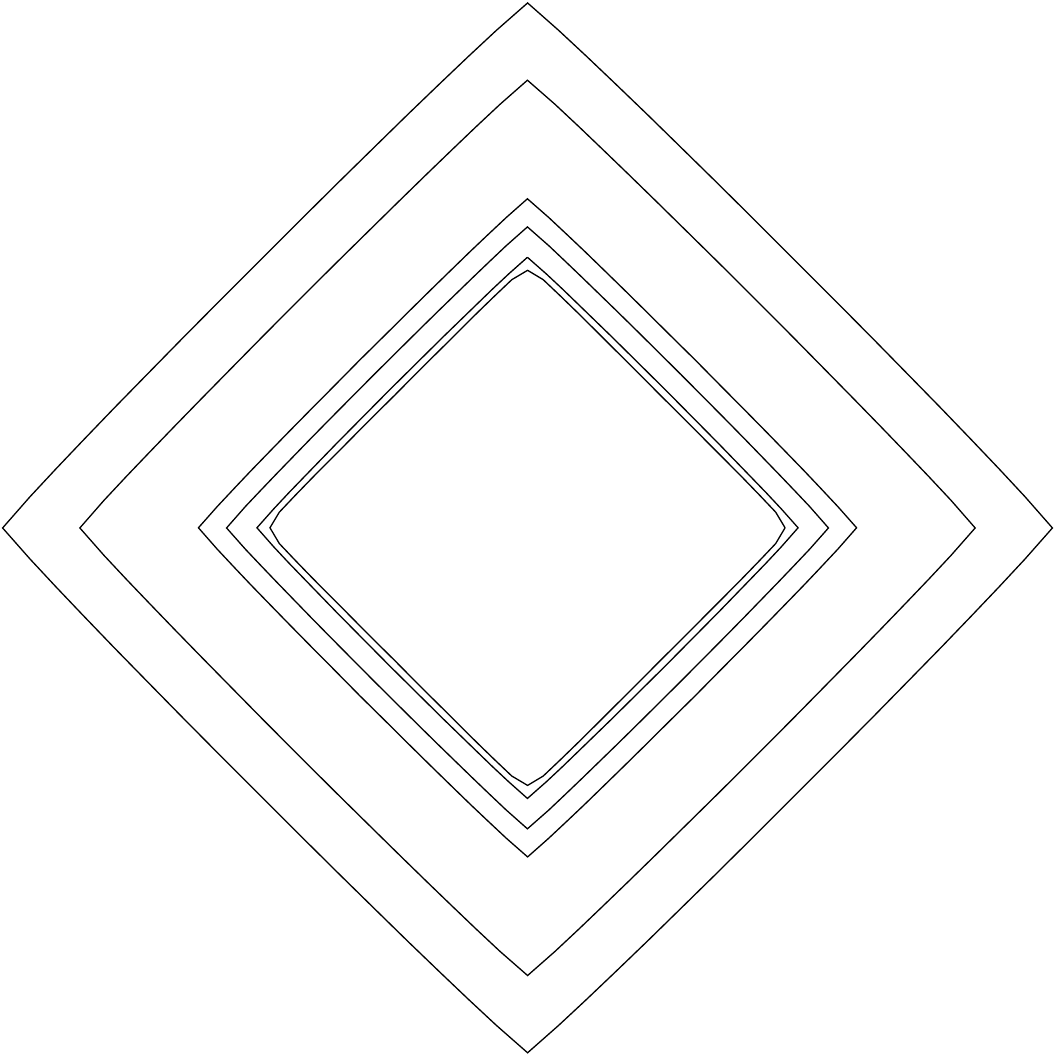} & \includegraphics[width=1\linewidth]{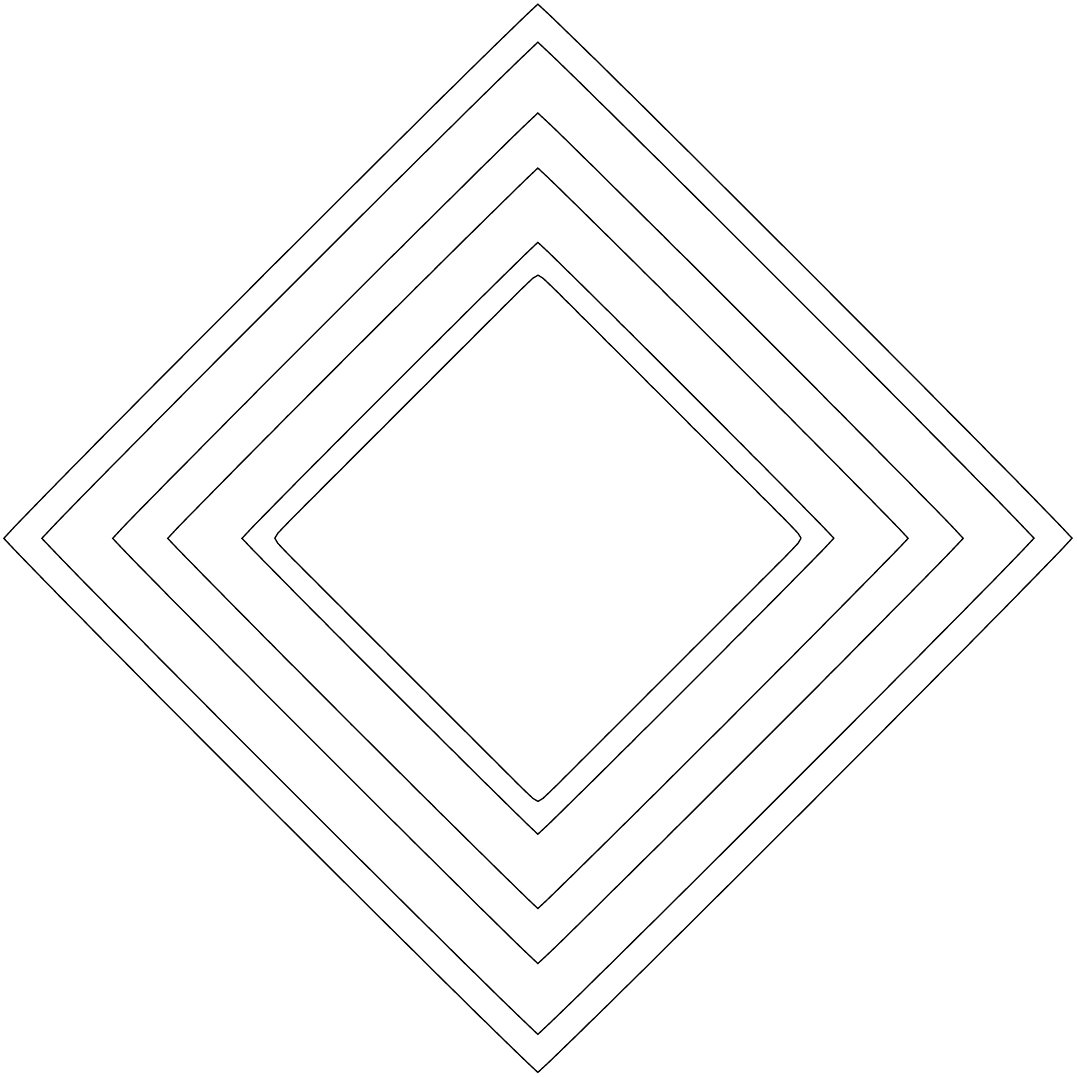} & \includegraphics[width=1\linewidth]{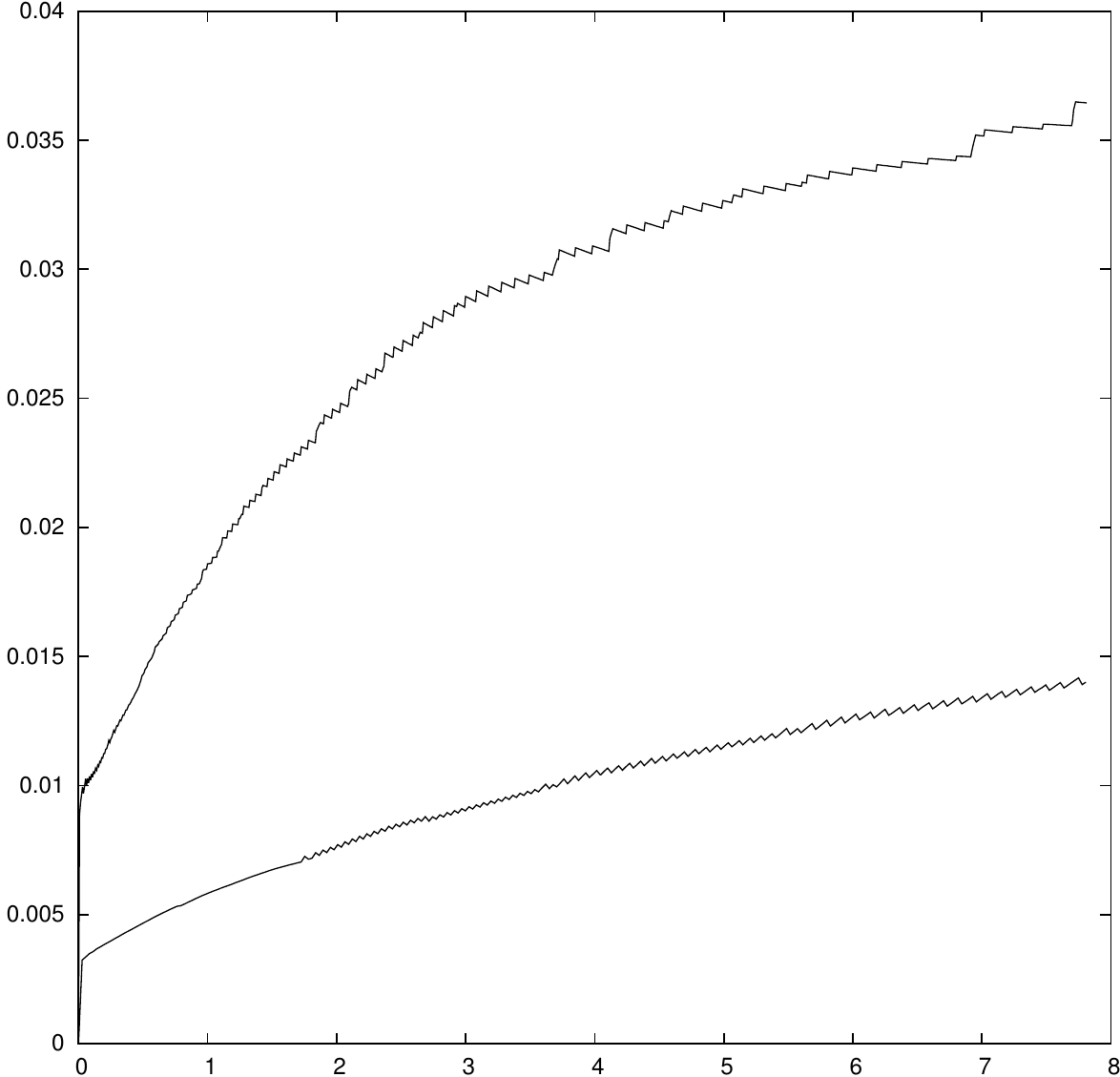}
\end{tabular}
\caption{\label{B8:fig:l1_evo} Evolution of the unit sphere with respect to the regularized $\ell^{1}$-norm under anisotropic Willmore flow for the anisotropy $\gamma(\cdot) = \|\cdot\|_{\ell_{\varepsilon}^{\infty}}$. The parameters are 
$h_0 = 0.0078$, $\varepsilon = 0.001$, $\tau = \tilde\tau =h_0^2$ 
and curves are plotted at times $t = 0, 10\tau, 50\tau, 100\tau, 500\tau, 1000\tau$ on the left 
and $h_0 = 0.0283$, $\varepsilon = 0.0001$, $\tau = h_0$, $\tilde \tau = h_0^2$, $t = 0, 10\tau, 50\tau, 100\tau, 200\tau, 275\tau$ in the middle.
On the right the associated $L^2$-errors are plotted over time, where the lower error curve corresponds to the evolution results on the left.
}
\end{figure}

\begin{table}[h]
\label{B8:tab:ellipse}
\begin{center}
\begin{tabular}{ | c || c | c | c || c | c | c |}
 \hline
  & & $L^2$-error & & & $L^2$-error & \\
  $n$   & $h(t)$ & $(\tau = \tilde{\tau} = h_0^2)$ & $eoc$ & $h(t)$ & $(\tau = \tilde{\tau} = h_0)$ & $eoc$  \\ \hline
  5     & 1.435e+0 & 1.648e-1 &       & 1.274e+0 & 1.942e-1 &  \\
  6     & 6.487e-1 & 3.476e-2 & 1.960 & 5.875e-1 & 7.089e-2 & 1.303 \\
  7     & 3.069e-1 & 8.762e-3 & 1.841 & 2.842e-1 & 3.424e-2 & 1.002 \\
  8     & 1.525e-1 & 2.182e-3 & 1.987 & 1.396e-1 & 1.724e-3 & 0.966 \\
\hline
\end{tabular}
\end{center}
\caption{\label{B8:table:eoc2} As in Table \ref{B8:table:eoc1} experimental orders of convergence are reported, now for the self-similar evolution of the ellipses 
(with half axis $6$ and $1$) under anisotropic Willmore flow. Here, again polygons with $2^n$ vertices are considered, equi-distributed along the initial ellipse with an initial grid size $h_0 = \frac{24.172}{2^n}$.
On the left the error is evaluated at time $t = 0.596576$ and on the right at time $t = 0.77238$.}
\end{table}
Now, we want to study crystalline anisotropies $\gamma(\cdot)= \|\cdot\|_{\ell^1}$ and 
$\gamma(\cdot)= \|\cdot\|_{\ell_\infty}$. As already pointed out, even though the formulation of the scheme itself doesn't explicitly need assumptions on the smoothness of $\gamma$, 
the application of the optimization algorithm requires the computation of derivatives 
of $\gamma$ up to order $3$. In fact, we use the following regularization:
For a small parameter $\varepsilon > 0$ we regularize the $\ell^1$-norm by
\begin{eqnarray*}
\ell^1_{\varepsilon}(z) = \sum\limits_{l =1}^2 \sqrt{\varepsilon |z|^2 + z_l^2}.
\end{eqnarray*}
Since in $\mathbb{R}^2$ the $\ell^{\infty}$-norm equals a rotated and scaled $\ell^{1}$-norm we use as regularization of the $\ell^{\infty}$-norm
\begin{eqnarray*}
\ell^{\infty}_{\varepsilon}(z) = \frac{\sqrt{\varepsilon |z|^2+ (z_1 + z_2)^2}}{2} + \frac{\sqrt{\varepsilon |z|^2 + (z_1 - z_2)^2}}{2}.
\end{eqnarray*}
Figure \ref{B8:fig:evogammalinf} shows the evolution of a sphere with respect to the regularized 
$\ell^{\infty}$-norm under the associated anisotropy Willmore flow with anisotropy
$\gamma(\cdot) = \|\cdot\|_{\ell^1_{\varepsilon}}$ for $\varepsilon = 0.0001$. 
Results on the self similar evolution of spheres with respect to the regularized 
$\ell^{1}$-norm  are depicted in Figure \ref{B8:fig:l1_evo}. In these simulations, we use the analog regularization for the dual anisotropy $\gamma^*$ required in the algorithm.

\begin{figure}
\begin{center}
\begin{tabular}{m{0.5\linewidth}m{0.2\linewidth}}
\includegraphics[width=1\linewidth]{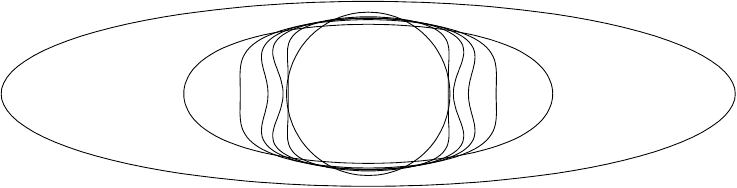} & \includegraphics[width=1\linewidth]{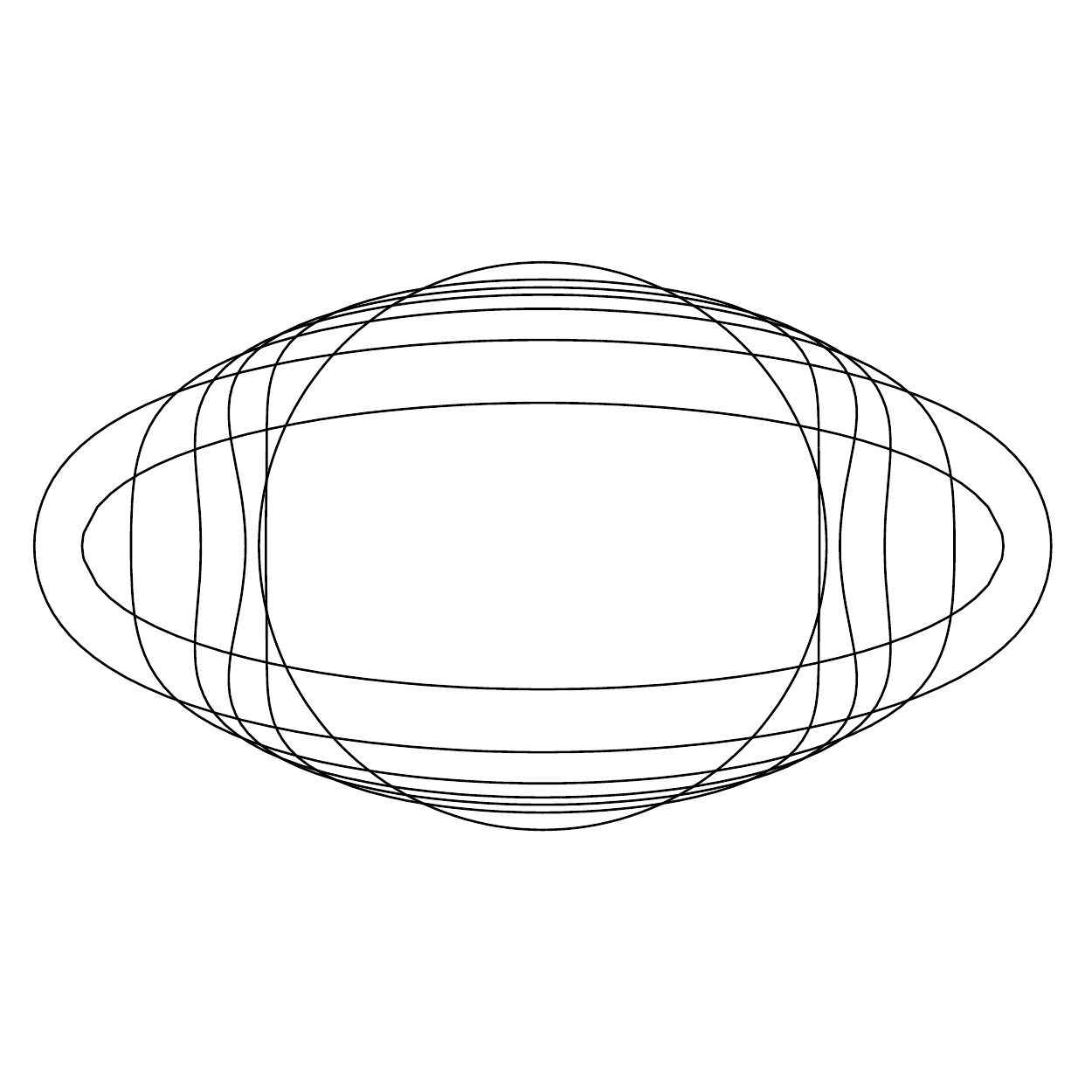} 
\end{tabular}
\caption{\label{B8:fig:lambda}
The impact of the parameter $\lambda$ is shown for the evolution of a circle to an ellipse with aspect ratio $4:1$ (i.e. $a_1 = 4$ and $a_2 = 1$). 
We evolve polygons with $160$ vertices approximating the unit sphere as initial curve, 
$h_0 = 0.0393$ and $\tau = \tilde\tau = 0.01$, $h = 0.000393$.
On the left $\lambda = 0.025$ and on the right $\lambda = 4$.}
\end{center}
\end{figure}

\begin{figure}
\begin{center}
\begin{tabular}{m{0.3\linewidth}m{0.3\linewidth}m{0.3\linewidth}}
\includegraphics[width=1\linewidth]{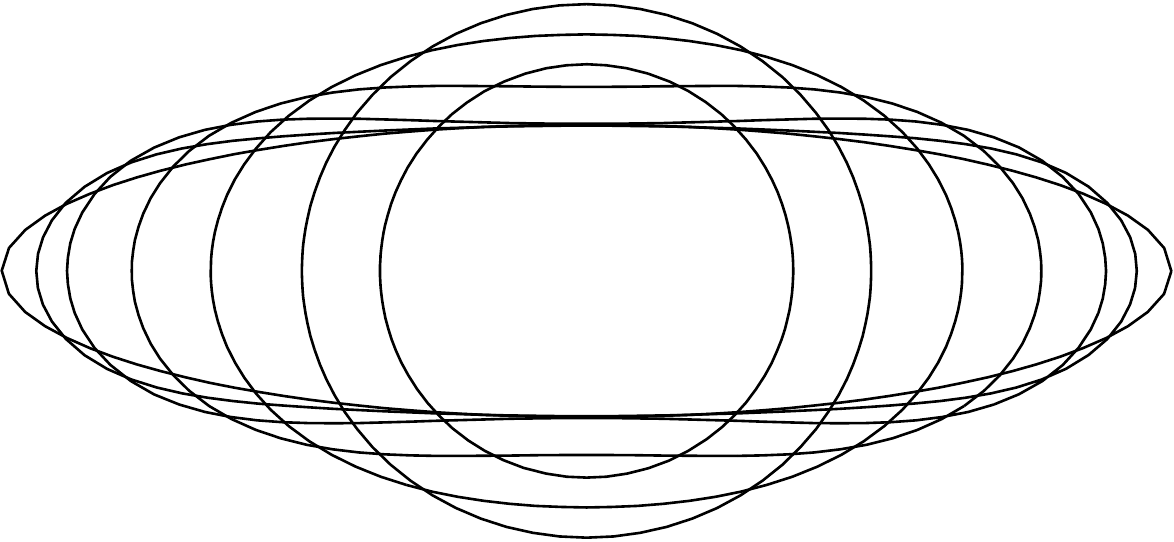} & \includegraphics[width=1\linewidth]{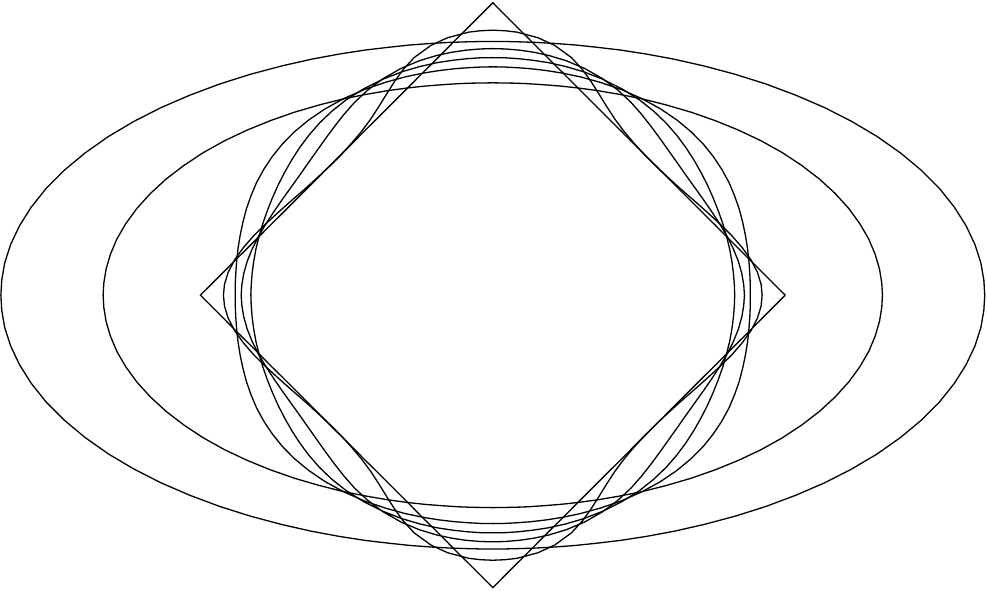} & \includegraphics[width=1\linewidth]{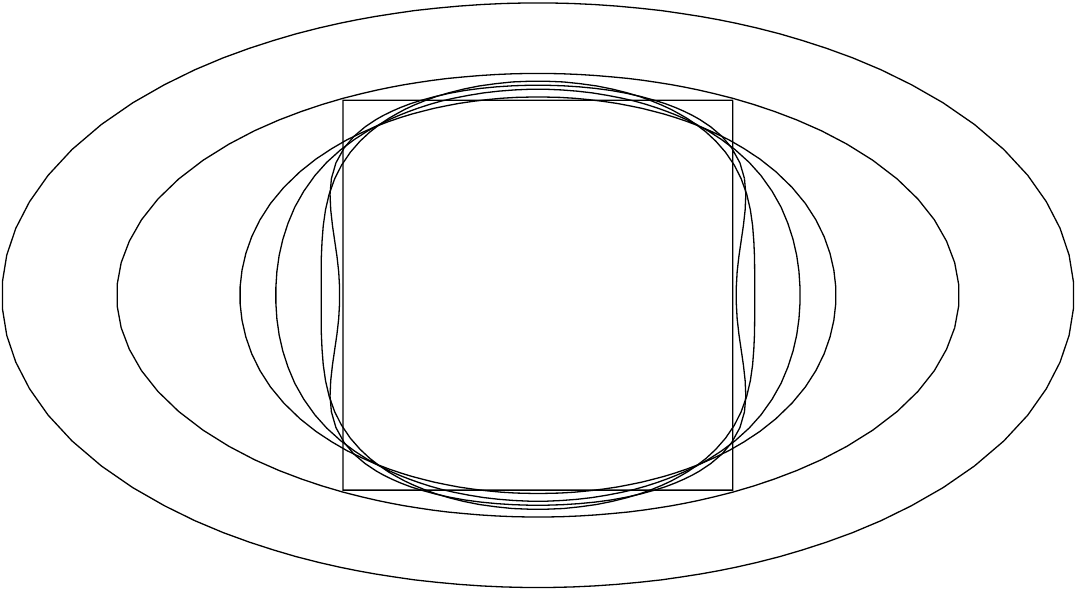}
\end{tabular}
\caption{\label{B8:fig:compare}The evolution of different initial shapes for different anisotropies is displayed. 
For all computations we use $100$ vertices and choose $\lambda = 0.25$. 
On the left we start with an ellipse with aspect ratio $4:1$ under an isotropic flow with $\gamma(\cdot) = \|\cdot \|$ ($h_0 = 0.1739$, $\tau = h_0$, $\tilde\tau = h_0^2$) results are shows at 
$t = 0, 0.1739, 0.5218, 1.739, 3.478, 6.956, 173.9$.
In the middle and on the right an ellipsoidal anisotropy with aspect ratio $2:1$ is used (i.e. $a_1 = 2$, $a_2 = 1$) 
in the first case (middle), we take as initial shape the unit sphere for the $l^1$-norm 
($h_0 = 0.0566$, $\tau = \tilde\tau = 0.001\,h_0$) and results are displayed at $t = 0, 0.00017, 0.00085, 0.00169, 0.006, 0.056, 0.251$).
In the second example (right), the initial shape is the unit sphere for the $l^{\infty}$-norm 
($h_0 = 0.08$ and $\tau = \tilde\tau= 0.01\,h_0$) and results are depicted for $t = 0, 0.0024, 0.008, 0.04, 0.08, 0.8, 4.8$.}
\end{center}
\end{figure}

Next, we generalize Willmore flow and replace the Willmore energy by
the modified energy
\begin{equation}
\label{B8:eq:WillmoreArea}
 e_\gamma[x]:= \int_{{\mathcal M}[x]} 
 \left(\frac{1}{2} \mathbf{h}_\gamma^2 + \lambda \right) \gamma(n[x]) d\mathit{a}\,,
\end{equation}
with a second term given by the anisotropic area weighted with a constant $\lambda>0$.
The incorporation of this generalized energy in our computational approach is straightforward.
The generalized flow combines expansive forcing with respect to the anisotropic Willmore flow of curves with contractive forcing due to the anisotropic mean curvature motion associated to the anisotropic area functional. Thus, for the generalized model we expect convergence 
to a limit shape given by a scaled Wulff shape, where the scaling depends on the factor 
$\lambda$. Fig. \ref{B8:fig:lambda} shows the impact of the factor $\lambda$ on the evolution,
whereas in Fig. \ref{B8:fig:compare} we compare the evolution of different initial shapes under the generalized anisotropic Willmore flow for different anisotropies.

\section*{Appendix}
Here, we collect the computational ingredients to evaluate the Lagrangian, its gradient and Hessian based on a standard local assembly procedure. 
In the following for vectors $x\in \mathbb{R}^{d+1}$  and functions $f$ we use the notation $f_{i,j}(x) = \frac{\partial f_i(x)}{\partial x_j}$ and in analogy for higher order derivatives. 
Furthermore, for matrices $A$ we use $f_{k,ij}(A) = \frac{\partial f_k(A)}{\partial A_{ij}}$ and again in analogy for higher order derivatives.
In fact, we can restrict ourselves to the local functionals
\begin{eqnarray}
{\mathbf{M}_{T,\gamma}}[Z,X] =  \frac{1}{(d+1)!} \left(\sum_{i =0,\ldots, d} \gamma^\ast(\bar Z_i)^2 \right) \gamma(R [\bar X])\,,~ 
{\mathbf{A}_{T,\gamma}}[X] =  \frac1{d!} \gamma(R [\bar X])\,,
\end{eqnarray}
where we denote by $\bar Z = (Z_0, \ldots, Z_{d})$ and $\bar X = (X_0, \ldots, X_{d})$ the corresponding vectors of simplex nodes in $\mathbb{R}^{d+1}$ with coordinate representation $Z_j = (Z_{jr})_{r=1,\ldots, d+1}$ and $X_j = (X_{jr})_{r=1,\ldots, d+1}$. 
Here, $R$ is a mapping from $\mathbb{R}^{(d+1)^2}$ to $\mathbb{R}^{d+1}$ representing the $90^\circ$ rotated edge vector for $d=1$ and the cross product of edge vectors for $d=2$, respectively.
For  $d=1$ we obtain for the first derivatives of $R[\bar X] = \PePoRuvectz{X_{02}-X_{12}}{X_{11}-X_{01}}$ 
with respect to the entries $(ij)$ with $i=0,\ldots,d$ and $j=1,\ldots, {d+1}$ 
$$
R_{,01}[\bar X] = \PePoRuvectz{0}{-1}\,,\quad R_{,02}[\bar X] = \PePoRuvectz{1}{0}\,,
\quad R_{,11}[\bar X] = \PePoRuvectz{0}{1}\,,\quad R_{,12}[\bar X] = \PePoRuvectz{-1}{0}\,.
$$
Because of the linearity of $R$ for $d=1$ all higher derivatives vanish.
For $d=2$ we have
\begin{eqnarray*}
R[\bar X] 
=\left( \sum_{u,v=1}^3 \epsilon_{iuv} (X_{1u}-X_{0u})(X_{2v}-X_{0v})\right)_{i=1,2,3}\,,
\end{eqnarray*} 
where $\epsilon_{wuv}$ is the Levi-Civita symbol ($\epsilon_{wuv}=\pm 1$ if $(w,u,v)$ is a even/odd permutation of $(1,2,3)$
and $0$ else). Thus, for $w=1,2,3$ we have
\begin{eqnarray*}
R_{w,js}[\bar X] &=&  \sum_{u,v=1}^3 \epsilon_{wuv}  \left(
(\delta_{1j}-\delta_{0j})\delta_{su}(X_{2v}-X_{0v}) + (\delta_{2j}-\delta_{0j})\delta_{sv} (X_{1u}-X_{0u})  \right),
 \\
R_{w,js\, lt}[\bar X] &=&  \sum_{u,v=1}^3 \epsilon_{wuv}  \left(
(\delta_{1j}-\delta_{0j})(\delta_{2l}-\delta_{0l})\delta_{su}\delta_{tv} + (\delta_{2j}-\delta_{0j}) (\delta_{1l}-\delta_{0l}) \delta_{sv} \delta_{tu} \right),
\\
\end{eqnarray*}
and all third derivatives $R_{u, ir \, js\, lt}[\bar X]$ vanish.
Here $j,l \in \{0,1,2\}$ refer to the local node and $s,t \in \{1,2,3\}$ to the spacial component.
Next we derive expressions for the derivatives of \eqref{B8:Mgamma} and \eqref{B8:Agamma} under the assumption that $\gamma, \gamma^{*}$ are sufficiently smooth  in $\mathbb{R}^{d+1} \setminus \{0 \}$ ( thus $Z_{i}, R[\bar X] \neq 0$):
\subsubsection*{Derivatives of ${\mathbf{M}_{T,\gamma}}$}
{\small \begin{eqnarray*}
\partial_{Z_{ir}} {\mathbf{M}_{T,\gamma}}[Z,X] &=& \frac{1}{(d+1)!} {\gamma}^{\ast,2}_{,r}(Z_{i}) \gamma(R [\bar X])\,, \\
\partial_{Z_{js}} \partial_{Z_{ir}} {\mathbf{M}_{T,\gamma}}[Z,X] &=& \frac{\delta_{ij}}{(d+1)!} 
{\gamma}^{\ast,2}_{,rs}(Z_{i}) \gamma(R [\bar X])\,, \\
\partial_{Z_{lt}} \partial_{Z_{js}} \partial_{Z_{ir}} {\mathbf{M}_{T,\gamma}}[Z,X] &=& \frac{\delta_{ij}\delta_{il}}{(d+1)!} {\gamma}^{\ast,2}_{,r s t}(Z_{i})\gamma(R [\bar X])\,, \\
\partial_{X_{ir}} {\mathbf{M}_{T,\gamma}}[Z,X] &=& \frac{1}{(d+1)!} \left(\sum_{\alpha =0,\ldots, d} {\gamma}^{\ast,2}(\bar Z_\alpha) \right)
\sum_{s=1}^m \gamma_{,s}(R [\bar X]) R_{s,ir}[\bar X]\,, \\[3ex] 
\partial_{X_{js}} \partial_{X_{ir}} {\mathbf{M}_{T,\gamma}}[Z,X] &=& \frac{1}{(d+1)!} \left(\sum_{\alpha =0,\ldots, d} {\gamma}^{\ast,2}(\bar Z_\alpha) \right) \\
&& \quad \cdot \left( \sum_{t=1}^m \gamma_{,t}(R [\bar X]) R_{t,ir\, js}[\bar X]  + \sum_{t,u=1}^m \gamma_{,t u}(R [\bar X]) R_{t,ir}[\bar X] R_{u,js}[\bar X]\right) \,, \\
\partial_{X_{js}} \partial_{Z_{ir}} {\mathbf{M}_{T,\gamma}}[Z,X] &=& \frac{1}{(d+1)!} {\gamma}^{\ast,2}_{,r}(Z_{i}) 
\sum_{t=1}^m \gamma_{,t}(R [\bar X]) R_{t,js}[\bar X]\,, \\
\partial_{X_{lt}} \partial_{Z_{js}} \partial_{Z_{ir}} {\mathbf{M}_{T,\gamma}}[Z,X] &=& \frac{\delta_{ij}}{(d+1)!} 
{\gamma}^{\ast,2}_{,r s}(Z_{i}) \sum_{u=1}^m \gamma_{,u}(R [\bar X]) R_{u,lt}[\bar X]\,, \\
\partial_{X_{lt}} \partial_{X_{js}} \partial_{Z_{ir}} {\mathbf{M}_{T,\gamma}}[Z,X] &=& \frac{1}{(d+1)!} 
 {\gamma}^{\ast,2}_{,r}(Z_{i}) \\
&& \quad \cdot  \left( \sum_{v=1}^m\gamma_{,v}(R [\bar X]) R_{v,js\, lt}[\bar X]
  + \sum_{v,u=1}^m \gamma_{,v u}(R [\bar X]) R_{v,js}[\bar X] R_{u,lt}[\bar X]\right)\,,
\end{eqnarray*}

\subsubsection*{Derivatives of ${\mathbf{A}_{T,\gamma}}$}
{\small
\begin{eqnarray*}
\partial_{Y_{ir}} {\mathbf{A}_{T,\gamma}}[X] &=& \frac1{d!} \sum_{s=1}^m \gamma_{,s}(R [\bar X]) R_{s,ir}[\bar X]
\,, \\
\partial_{Y_{js}} \partial_{Y_{ir}} {\mathbf{A}_{T,\gamma}}[X] &=& \frac1{d!}  \left(\sum_{t=1}^m \gamma_{,t}(R [\bar X]) R_{t,ir\, js}[\bar X]  + \sum_{t,u=1}^m \gamma_{,t u}(R [\bar X]) R_{t,ir}[\bar X] R_{u,js}[\bar X]\right)\,, \\
\partial_{Y_{lt}} \partial_{Y_{js}} \partial_{Y_{ir}} {\mathbf{A}_{T,\gamma}}[X] &=& \frac1{d!} \Big( \sum_{u,v=1}^m
\gamma_{,v u}(R [\bar X]) R_{v,ir\, js}[\bar X] R_{u,lt}[\bar X]  + 
 \sum_{v=1}^m\gamma_{,v}(R [\bar X]) R_{v,ir\, js\, lt}[\bar X]  \\[-.8ex]
&&   +  \sum_{u,v,w=1}^m\gamma_{,v u w}(R [\bar X]) R_{v,ir}[\bar X] R_{u,js}[\bar X] R_{w,lt}[\bar X]\\[-.8ex]
&&  +  \sum_{u,v=1}^m\gamma_{,v u}(R [\bar X]) R_{v,ir\, lt}[\bar X] 
  R_{u,js}[\bar X] 
    +  \sum_{u,v=1}^m\gamma_{,v u}(R [\bar X]) R_{v,ir}[\bar X] R_{u,js\, lt}[\bar X]
\Big) \,. \\
\end{eqnarray*}
}

\subsection*{Acknowledgement.}
Ricardo Perl was supported by the DFG project Ru 567/14-1 and 
Martin Rumpf acknowledges support by the SFB 611.

\bibliographystyle{abbrv}
\bibliography{../Bibtex/own,../Bibtex/all,../paola}
\end{document}